\documentclass[journal]{IEEEtran}
\pdfoutput=1                                                     
\usepackage{amsmath,amsfonts,amssymb,mathtools,bm,cite}
\usepackage[normalem]{ulem}
\usepackage{graphicx}

\usepackage{color}
\newcommand{\be}{\begin{equation}}
\newcommand{\ee}{\end{equation}}

\IEEEoverridecommandlockouts                              

\title{\LARGE \bf Fuel Minimisation for a Vehicle Equipped with a Flywheel and Battery on a Three-Dimensional Track} 
\author{Mehdi Imani Masouleh and David~J.~N.~Limebeer
\thanks{This work was supported by the UK Engineering and Physical Sciences Research Council.}
\thanks{Mehdi Imani Masouleh (mehdi.imanimasouleh@eng.ox.ac.uk) and David~J.~N.~Limebeer (david.limebeer@eng.ox.ac.uk) are with the Department of Engineering Science, University of Oxford, Parks Road, Oxford, OX1 3PJ.}%
}

\begin{document}
\maketitle
\thispagestyle{empty}
\pagestyle{empty}

\begin{abstract}
An optimal control based methodology is proposed for minimising the combustible fuel consumption of a hybrid vehicle equipped with an internal combustion engine, a high-speed flywheel and a battery.
The three-dimensionality of the road is recognised by the optimal control calculations.
Fuel efficiency is achieved by optimally exploiting the primary and secondary energy sources and controlling the vehicle so that the fuel consumption is minimised for a given, but arbitrary three-dimensional route. A time-of-arrival constraint rather than a driving cycle is  used. The benefits of using multiple auxiliary storage systems are demonstrated and a lower-bound estimate of the fuel consumption is presented.
\end{abstract}


\section{Introduction} 
The low efficiency and the air-pollution side effects associated with internal combustion engine (ICE) usage are well known. As is now widely appreciated, several of the problems endemic to ICEs can be mitigated using a variety of secondary energy 
storage technologies such as fuel cells, lithium-based battery systems, high-speed flywheels and supercapacitors. In broad terms, fuel cells and lithium-based battery systems have good energy storage properties, while high-speed flywheels and supercapacitors can be utilised for their high power density characteristics \cite{Doucette2011,Guzzella2007}. As a result of their poor power density properties, battery packs in commercial electric vehicles tend to be oversized. Other battery-related issues include long charging times and shortened life expectancy, especially when they are cycled at high charge and discharge rates. Supercapacitors and flywheel-based systems are an ideal complement to batteries, because of their high power density characteristics under both charging and discharging. The fundamentals and applications of fuel cells, including the main reactions, are reviewed in \cite{Carrette2000}. A comprehensive review of lithium-ion battery technologies is given in \cite{Horiba2014}. The physical structure of some flywheel-based systems is reviewed in a vehicle context in \cite{Dhand2013}. Hybrid energy storage systems (HESS) based on batteries and supercapacitors are reviewed in \cite{Ju2014} and the references therein.

Aside from the technologies themselves, energy storage modelling as well as their optimal deployment are important issues.  The modelling and control of hybrid electric vehicles (HEVs) are reviewed in \cite{Shen2011}, where various powertrain topologies and control strategies are discussed. It is pointed out that global optimization can be used as a `what's possible' benchmark  for evaluating energy management strategies. Another good survey paper reviewing 180 papers on the optimal energy management of HEVs and plug-in HEVs can be found in \cite{Panday2014}. Power management control strategies can be divided into offline and online methodologies. Online management strategies include methods such as look-up tables, state machines, thermostat control, {\em Equivalent Consumption Minimisation Strategies} (ECMS) \cite{Paganelli2002,Rousseau2008,Sciarretta2004}, Neural networks \cite{Baumann1998}, particle swarm optimisation \cite{Wang2006}, model predictive control (MPC) \cite{West2003} and fuzzy control \cite{Lee1998}. Offline strategies that often focus on global optimisation include dynamic programming (DP) \cite{Lin2003,Rousseau2008}, linear programming \cite{Tate2000}, nonlinear programming \cite{Perez2009}, Stochastic DP  \cite{liu2008, Moura2011, Kim2007} and genetic algorithms \cite{Piccolo2001}. An evolutionary algorithm was applied in \cite{Qi2016} over a sliding window to allow real-time power management and a data-driven reinforcement learning algorithm was proposed in \cite{Qi2016data}. ECMS does not guarantee charge sustainability and hence an Adaptive ECMS (A-ECMS) algorithm was introduced in \cite{Musardo2005} to update the equivalence factor `on-the-fly' on the basis of past and predicted driving conditions.

The majority of the work in the literature uses driving cycles as benchmarks for performance evaluation. As was recognised in \cite{Roy2014}, optimising vehicle parameters over one drive cycle does not necessarily mean that the vehicle will perform well on other drive cycles. In \cite{Roy2014}, a method was proposed for optimising an HEV over a range of drive cycles with different levels of driving aggressiveness and traffic conditions, in order to reduce the fuel economy variability with respect to drive cycle changes. While seeking to address the limitations associated with single driving cycle usage, this method is still restricted by the driving cycle combination used and there were cases where a single driving cycle resulted in a lower fuel variability compared with the proposed multi-cycle method.
 
The (combustible) fuel used by any vehicle will depend on the vehicle's speed, with higher average speeds typically resulting in a higher fuel usage. A contribution of this work is proposing a method of minimising fuel consumption over different driving conditions without using drive cycles. The idea is to simultaneously optimise the powertrain’s energy deployment and driving strategy over a given (but arbitrary) route. Key in this procedure is the selection of a time-of-arrival constraint, which acts as an `aggressiveness' surrogate. A short travel time corresponds to aggressive driving and a generally higher fuel consumption. An optimal control algorithm then seeks to minimise the fuel consumption, while ensuring a `just-in-time' arrival. The optimal control calculation makes use of a realistic vehicle model, with three degrees of freedom, and a non-linear tire model. It is shown that flywheels are an excellent means of reducing fuel consumption in manoeuvres where high levels of braking are involved, such as extra-urban driving on fast roads.

Another distinguishing feature of this work over the majority of the existing literature is that unlike rule-based methods, system dynamics and non-linear constraints can be included explicitly as part of the optimal control problem and optimal rather than near optimal results are obtained. A pseudo-spectral method will be used to solve the optimal control problem, which makes use of first- and second-order gradient information for fast convergence. Compared to DP/SDP, models of much higher complexity can be solved using this approach.

In this paper, we will take the road’s three-dimensional geometric profile into consideration. In mechanical systems with comparable potential and kinetic energies, the optimal control strategy is strongly dependent on trade-offs between the two energy sources, as was demonstrated by the famous `Minimum Time to Climb' problem associated with jet-powered aircraft \cite{Kaiser1944}. In Section\,\ref{IllustrativeExample} a simple motivating example is given that highlights the importance of changes in the road gradient. In Section\,\ref{CarAndTrackModel} the vehicle and track models employed in this study are described. In Section\,\ref{OptimalControl} the optimal control problem is cast in a standard form and the numerical method used to solve it is described. Numerical results are presented and discussed in Section\,\ref{Results}. The conclusions are given in Section\,\ref{Conclusion} and the vehicle parameters used in the study are provided in the Appendix. 

\section{An Illustrative Example}\label{IllustrativeExample}
Imagine a path-constrained point mass (bead) that is required to reach a given destination within a pre-specified time. 
Suppose the trajectory's starting point is $(0,0)$ in absolute Cartesian coordinates and that the path is parabolic and defined by $y=ax^2+bx$.
If the destination point $(x_f,\,\star)$ is given, there holds $y_f=ax_f^2+b x_f$ and so the path is specified by a single free parameter. 
We would like to minimise the energy supplied to the bead in order that it arrives at the destination `just in time'.
If the effects of air resistance and friction forces are neglected, the bead's height is indicative of its speed.
If the bead has unity mass, conservation of energy dictates
\begin{equation}
\frac{1}{2}v^2=E_{fuel}-gy.
\label{energyConservationEqn}
\end{equation}
The bead's speed is $v$, its instantaneous height above the origin is $y$ and $E_{fuel}$ is the amount of external energy supplied. The speed of the bead is thus
\begin{equation}
v = \sqrt{2E_{fuel}-2gy}.  \label{vsqrt}
\end{equation}
 
An infinitesimal path segment can be described as
\begin{align}
ds &= \sqrt{dx^2+dy^2} = dx \sqrt{(1+\left( \frac{dy}{dx} \right)^2} \nonumber \\
   &=dx \sqrt{1+(2ax+b)^2}. \label{infinitesimalLength}
\end{align}

Using \eqref{vsqrt} and \eqref{infinitesimalLength} the manoeuvre time is given by 
\begin{align}
T &=\int\limits_{0}^{T} dt = \int\limits_{0}^{L} \frac{1}{v}ds \\
 &= \int\limits_{0}^{x_f} \sqrt{\frac{1+(2ax+b)^2} {2E_{fuel}-2g(ax^2+bx)}} \,\, dx.
 \label{timeEquation}
\end{align} 

The thrust programme that minimises the total external energy supplied is the optimal control problem of minimising
\begin{align}
E_{fuel}=\int\limits_{s_0}^{s_f} F(s) ds = \int\limits_{0}^{t_f} F(t) v(t) dt
\end{align}
subject to state dynamics
\begin{eqnarray}
  \left\{
  \begin{aligned}
	\dot{x}& = v \frac{dx}{ds} = \frac{v}{\sqrt{1+\left( \frac{dy}{dx} \right)^2} } \\
	\dot{v}& =F - \frac{g \frac{dy}{dx}}{\sqrt{1+\left( \frac{dy}{dx} \right)^2} }
  \end{aligned}
  \right.
\end{eqnarray}
where $F$ is the externally applied force. The boundary conditions are
\begin{eqnarray}
  \left\{
  \begin{aligned}
	x(t_0)&=0 \quad x(t_f)=5 \\
	v(t_0)&=0 \quad v(t_f)= \mbox{free}.	
  \end{aligned}
  \right.
\end{eqnarray}
The control Hamiltonian of this system is
\begin{align}
\mathcal{H}=\lambda_1 \dot{x} + \lambda_2 \dot{v} + Fv.
\end{align}

Since $v(t_f)$ is free, $\lambda_2(t_f)=0$. For the quadratic path the control Hamiltonian becomes
\begin{align}
\mathcal{H}=F(\lambda_2+v) + \frac{\lambda_1v}{\sqrt{1+ (2ax+b)^2}} + \frac{\lambda_2 g(2ax+b)}{\sqrt{1+ (2ax+b)^2}} 
\label{hamiltonian}
\end{align}
with co-state equations
\begin{align}
\dot{\mathbf{\lambda}}=-\frac{\partial \mathcal{H}}{\partial \mathbf{x} } = \begin{pmatrix}  \partial \mathcal{H}/\partial x\\-F -\lambda_1/\sqrt{1+(2ax+b)^2}\end{pmatrix}
\end{align}
where
\begin{align}
\partial \mathcal{H}/\partial x&=\frac{2\lambda_2 ag}{\sqrt{1+ (2ax+b)^2}} -\frac{2\lambda_2 ag(2ax+b)^2}{(1+ (2ax+b)^2)^{3/2}} \nonumber\\
& + \frac{2\lambda_1 av(2ax+b)}{(1+ (2ax+b)^2)^{3/2}}.
\end{align}

\begin{figure}[tb]
\centering
		\includegraphics[width=0.5\textwidth]{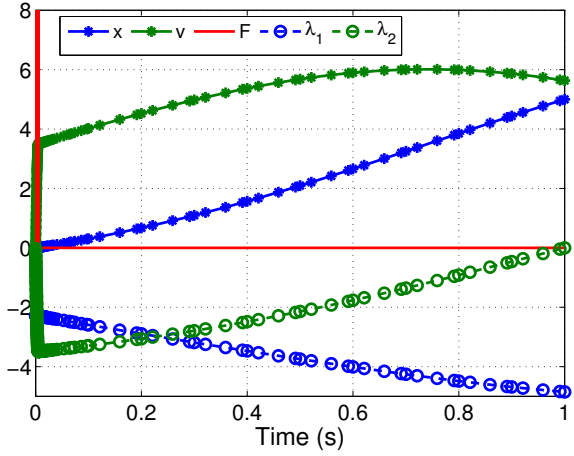}
\caption{Solution of the optimal control problem for $a=0.1$ and $y_f=-1$.}
\label{ocSol}
\end{figure}

Since the drive force is bounded, Pontryagin's minimum principle determines that
\begin{eqnarray}
  \left\{
  \begin{aligned}
	F=F_{max} \quad \text{if}  \quad \lambda_2 < -v \\
	F=0 \quad \text{if}  \quad \lambda_2 > -v\\
	F=singular \quad \text{if}  \quad \lambda_2 = -v.
  \end{aligned}
  \right.
\end{eqnarray}
In a manner reminiscent of simple variants of the Goddard rocket problem, the optimal thrust programme turns out to be impulsive \cite{Tsien1951}, with all the external energy injected into the system at the beginning of the manoeuvre. To illustrate this the described optimal control problem was solved with GPOPS in the case that
$y_f=-1$ and $a=0.1$; the states, control and co-states are shown in Fig.\,\ref{ocSol}. 

\begin{figure}[ht]
\centering
\includegraphics[width=0.47\textwidth]{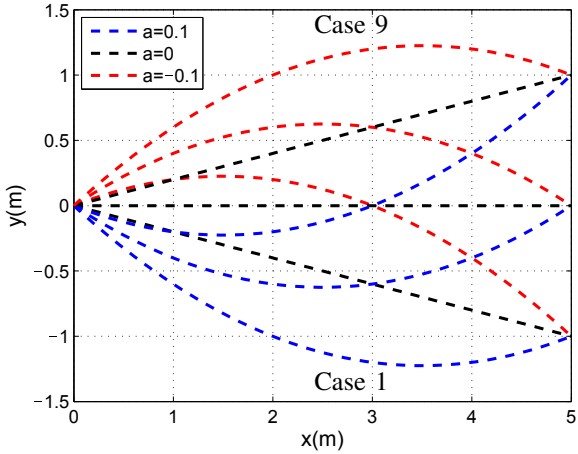}
\put(-120,+25){Case 1}
\put(-120,+160){Case 9}
\caption{Quadratic paths taken by the bead from the starting point $(0,0)$, to the terminal points $(5,-1)$, $(5,0)$ and $(5,1)$. }
\label{quadraticPaths}
\end{figure}

Nine cases are considered that highlight the significance of interactions between the kinetic and potential energy of the bead. These paths correspond to three values for the $a$ parameter and three values for the terminal point; the resulting paths are shown in Fig.\,\ref{quadraticPaths}.

For each of the nine cases, Table\,\ref{pathParams} shows the path lengths $L$, 
the zero-fuel times $T_0$ (in the cases that a zero-fuel journey is possible), and the external energy $E_{fuel}$ required to complete the journey under a 1\,s time-of-arrival constraint.
\begin{table}[ht]
\begin{center}
\caption{Path lengths ($L$), zero-fuel arrival times ($T_0$) and fuel requirements ($E_{fuel}$) for a 1\,s arrival time.}
\begin{tabular}{|c|c|c|c|c|c|c|}
\hline
Case & $y_f$ & $a$ & $b$ &  $L (m)$  & $T_0 (s) $ & $E_{fuel} (J)$ \\ 
\hline
1&-1 & 0.1 &  -0.7  &  5.29 & 1.57 	& 6.02 \\ \hline
2&-1 & 0   &  -0.2  &  5.10 & 2.30	&  8.56 \\ \hline
3&-1 & -0.1&  0.3   &  5.29 & -		& 13.51 \\ \hline 
4& 0 & 0.1 &  -0.5  &  5.20 & 2.38  & 9.70 \\ \hline
5& 0 & 0   &  	0   &  5.00    &   &  12.50 \\ \hline
6& 0 & -0.1 &  0.5  &  5.20    & -  & 17.73 \\ \hline 
7&+1 & 0.1 &  -0.3  &  5.29    & -  & 15.83 \\ \hline
8&+1 & 0   &  -0.2  &  5.10    & -  &  18.37 \\ \hline
9&+1 & -0.1 &  0.7  &  5.29    & -  & 23.31 \\ \hline 
\end{tabular}
\end{center}
\label{pathParams}
\end{table}

The second column in Table\,\ref{pathParams} shows the elevation of the terminal point. To calculate the arrival time, \eqref{timeEquation} is solved numerically by running a bisection on $E_{fuel}$ until the journey time constraint is met.

Superior fuel consumption performance is achieved when the initial part of the journey is downhill, because the early conversion of potential energy into kinetic energy results in higher speeds throughout the journey. Shorter journey paths do not necessarily result in lower fuel consumption (e.g. compare Case 1 and 2). As one would expect, the journeys in Cases\,7 to 9 are the most arduous interms of fuel consumption, due to elevated terminal points. This example demonstrates that the path elevation curvature can have a significant influence on fuel usage and sometimes in a counter-intuitive manner.

\vspace{+5mm}
\section{The Mathematical Model} \label{CarAndTrackModel}
\vspace{+5mm}
A bicycle model of the car is used that has yaw, lateral and longitudinal freedoms, and nonlinear tires. The road is assumed three-dimensional and is represented by a geometric construct called a `ribbon', which is describable in terms of three curvature variables \cite{Perantoni2015Track}.

\subsection{Track Model} \label{Track_Model}
\begin{figure}[ht]
\centering
\includegraphics[width=0.5\textwidth]{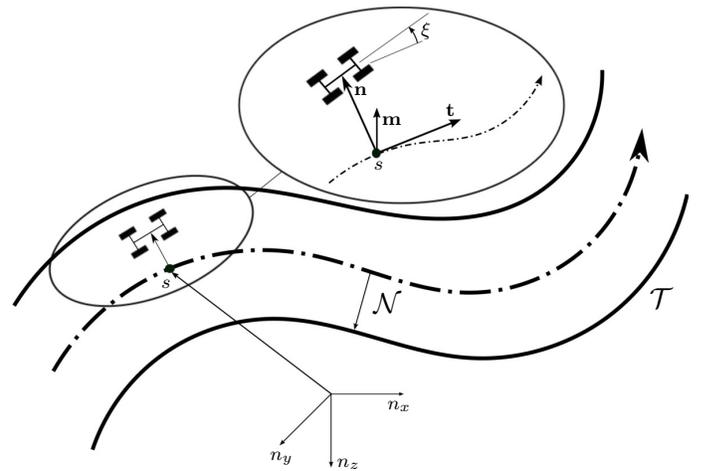}
\caption{Differential-geometric description of a track segment $\cal T$. The independent variable $s$ represents the distance travelled from the start line. The track half-width is ${\cal N}$, with $\xi$ the car's yaw angle relative to the track spine's tangent direction. The inertial reference frame is given by $n_x$, $n_y$ and $n_z$.}
\label{path_plan}
\end{figure}

A moving coordinate system (called a Darboux frame) is used to describe the track. As shown in Fig.\,\ref{path_plan} the origin of this moving system is `dragged along' by the car. The independent variable in the track description is $s$, which is the distance travelled by the car from some starting point, projected onto the track spine. The spine could be, but need not be, the track centre line.

The Darboux frame is described by the orthogonal moving triad $[{\bm t}\,{\bm n}\,{\bm m}]$, with ${\bm m}$ normal to the road surface.
The road is represented locally by the $\bm{t}$-$\bm{n}$ plane, which `travels' with the car down the road.

The orientation of this local `patch' is described in terms of three Euler angles that are all functions of $s$ \cite{Perantoni2015Track}. If the track orientation is described by the roll, yaw and pitch angles, respectively $\phi$, $\mu$ and $\theta$, then the track curvatures are given by 
\begin{equation}
     \bm{\Omega}=\begin{bmatrix}
	    \Omega_x \\
	    \Omega_y \\
	    \Omega_z
    \end{bmatrix}
    =\frac{1}{\dot{s}} \begin{bmatrix}
    	    \dot{\phi} - \sin(\mu) \dot{\theta}  \\
    	    \cos(\phi) \dot{\mu} +\cos(\mu)\sin(\phi) \dot{\theta} \\
    	     -\sin(\phi) \dot{\mu} + \cos(\mu) \cos(\phi) \dot{\theta} 
        \end{bmatrix}.
\end{equation}
The angular velocity of the Darboux frame is given by
\begin{equation}
     \bm{\omega}= \begin{bmatrix}
	    \omega_x &   \omega_y &   \omega_z
    \end{bmatrix}^T
    =\dot{s}\bm{\Omega}.
\end{equation}

The next kinematic relationships we will require relates to the way in which the car progresses down the road.
Suppose that the absolute velocity of the car in its body-fixed coordinate system is $[u~v~w]^T$; the longitudinal velocity component $u$ is determined by the throttle/brakes, the lateral component $v$ is determined by the steering, while the vertical component $w$ is determined by the track characteristics. The car's geometric centre (the car's mass centre projected down on to the road) is given by ${\boldsymbol n}=[0~n~0]^T$ in the Darboux frame.
The absolute velocity of the car in the Darboux frame is:
\begin{align}
\begin{bmatrix}
\dot{s}\\
\dot{n}	\\
0 
\end{bmatrix}
={\boldsymbol n} \times {\boldsymbol \omega} + R_z(\xi){\boldsymbol v} =\begin{bmatrix}
n \omega_z + u \cos \xi - v \sin \xi	\\
v \cos \xi + u \cos \xi	\\
	w - n  \omega_x
\end{bmatrix}
\label{kinematics1}
\end{align}
where
\begin{equation}
R_z(\xi) = R(\textbf{e}_z, \xi)=\left[
\begin{array}{rrr}
	\cos \xi & -\sin \xi & 0 \\
	\sin \xi  & \cos \xi & 0 \\
	0 & 0 & 1 \\
\end{array} \label{caryaw}
\right]
\end{equation}
represents the yawing of the car relative to the spine of the track. 
The first term in (\ref{kinematics1}) derives from the angular velocity of the Darboux frame, while the second is the velocity of the car's geometric centre expressed in Darboux coordinates. The first row of (\ref{kinematics1}) gives the speed of the origin of the Darboux frame in its tangent direction and can be re-written as
\begin{equation}
\dot{s} = \frac{u \cos \xi - v \sin \xi}{1 - n \Omega_z}, \label{speed_s}
\end{equation}
since $\omega_z=\dot{s} \Omega_z$. The function
\begin{equation}
S_f (s) = \frac{dt}{ds}
\end{equation}
transforms `time' as the independent variable into the `elapsed distance' as the independent variable; it is assumed that  $S_f (s)$ and its inverse are non-zero everywhere on the spine curve.
The second row of (\ref{kinematics1}) is
\begin{equation}
\dot{n} = u \sin \xi + v \cos \xi,  \label{n_rate}
\end{equation}
which describes the way the vehicle moves normal to the spine. Transforming \eqref{n_rate} into the distance-travelled domain gives
\be
n^\prime = S_f (s) \left( u \sin \xi  + v \cos \xi  \right)  \label{n_rate_s}
\ee
in which $n^\prime$ is the derivative of $n$ with respect to $s$.
Expressing the absolute angular velocity of the car in its body-fixed reference frame gives
\be
\bar{\boldsymbol \omega} = \left[
\begin{array}{c}
\bar{\omega}_{x} \\
\bar{\omega}_{y} \\
\bar{\omega}_{z}
\end{array}
\right] = 
\left[
\begin{array}{c}
\cos \xi \omega_x + \sin \xi \omega_y \\
\cos \xi \omega_y - \sin \xi \omega_x \\
\omega_z + \dot{\xi} \label{carabsz}
\end{array}
\right].
\ee
This will be used in the next section to derive the vehicle's equations of motion. The car's yaw angle in Darboux frame $\xi$ expressed in distance domain is deduced from the third row of (\ref{carabsz}) by integrating
\begin{eqnarray}
\xi^\prime = S_f (s) \bar{\omega}_{z} -\Omega_z. 
\label{xi_rate_s}
\end{eqnarray}
\subsection{Dynamics} \label{dynamics}
\begin{figure}[tb]
	\centering
	\includegraphics[width=0.47\textwidth]{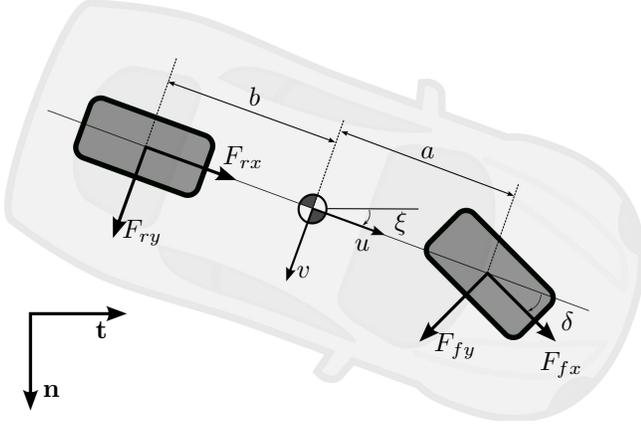}
    \put(-85,100){$a$}%
	\put(-150,120){$b$}%
    \put(-208,32){$\bf t$}%
    \put(-228,10){$\bf n$}%
    \put(-110,65){$u$}%
    \put(-132,55){$v$}%
    \put(-95,73){$\xi$}%
    \put(-32,35){$\delta$}%
    \put(-40,20){$F_{fx}$}%
    \put(-80,27){$F_{fy}$}%
    \put(-160,98){$F_{rx}$}%
    \put(-198,70){$F_{ry}$}%
\caption{Tire force system. The car's yaw angle with respect to the Darboux frame, which is defined in terms of the vectors $\bf t$ and $\bf n$, is $\xi$ and $\delta$ is the steering angle.}
\label{tyre_force_diag}
\end{figure}

The equations describing the dynamics of the car are derived using standard vectorial methods. 
The absolute velocity of the car's mass centre (expressed on the vehicle's coordinate system) can be written as
\begin{align}
{\boldsymbol v}_{B} &= {\boldsymbol v} + \bar{\boldsymbol \omega} \times {\boldsymbol h} \nonumber
\end{align}
\begin{align}
&= \begin{bmatrix}
u\\
v	\\
n \omega_x
\end{bmatrix} + 
\begin{bmatrix}
\bar{\omega}_x \\
\bar{\omega}_y \\
\bar{\omega}_z
\end{bmatrix} \times \begin{bmatrix}
0 \\
0 \\
- h \end{bmatrix}
= \begin{bmatrix}
u - h \bar{\omega}_y \\
v + h \bar{\omega}_x \\
n \omega_x
\end{bmatrix}
\end{align}
where $\times$ denotes the cross product. The Newton-Euler equations for this system are given by
\begin{eqnarray}
M \left(\dot{\boldsymbol v}_B + \bar{\boldsymbol \omega} \times {\boldsymbol v}_B \right) &=& \bm{F}_B + MgR^T {\boldsymbol e}_z \label{Fbal} \\
I_B \dot{\bar{\boldsymbol \omega}} + \bar{\boldsymbol \omega} \times (I_B \bar{\boldsymbol \omega}) &=& \bm{M}_B, \label{Mbal}
\end{eqnarray}
where the car's inertia matrix is assumed to be diagonal and is given by $I_B = \mbox{diag} (I_x~I_y~I_z )$, with $\bm{F}_B = [F_x~F_y~F_z ]^T$ and $\bm{M}_B = [M_x~M_y~M_z ]^T$ being the external force and moment. The last term in \eqref{Fbal} is due to the gravitational acceleration of the car's mass centre and can be expressed as
\begin{align}
MgR^T {\boldsymbol e}_z & = R^T_z(\xi) R^T_x(\phi) R^T_y(\mu) \begin{bmatrix}
0 & 0 & Mg 
\end{bmatrix}^T \nonumber \\
&= Mg \begin{bmatrix}
\sin \xi \sin \phi \cos \mu - \cos \xi \sin \mu \\
\sin \xi \sin \mu + \cos \xi \sin \phi \cos \mu \\
\cos \phi \cos \mu
\end{bmatrix}.
\label{Fgrav}
\end{align}
The car's equations of motion can now be assembled from (\ref{Fbal}), (\ref{Mbal}) and (\ref{Fgrav}) as follows:
\begin{align}
\dot{u} &= (v+h \bar{\omega}_x) \bar{\omega}_z - n \omega_x \bar{\omega}_y \nonumber + h \dot{\bar{\omega}}_y \\
& \quad + g \left(\sin \xi \sin \phi \cos \mu - \cos \xi \sin \mu \right) + F_x/M \label{longF}  \\
\dot{v} &= n \omega_x \bar{\omega}_x - (u-h \bar{\omega}_y) \bar{\omega}_z  - h \dot{\bar{\omega}}_x \nonumber  \\
& \quad + g \left( \sin \xi \sin \mu + \cos \xi \sin \phi \cos \mu \right) + F_y/M \label{latF} \\
\dot{\bar{\omega}}_z &= \left( (I_x-I_y)\bar{\omega}_x \bar{\omega}_y  + M_z \right)/I_z,  \label{yaw_bal}
\end{align}
in which $F_x$, $F_y$ are the resultant longitudinal and lateral forces, and $M_z$ is the z-axis tire moment acting on the car. These quantities are given by
\begin{align}
 F_x &= F_{fx}\cos \delta - F_{fy}\sin \delta  + F_{rx} + F_{ax}  \nonumber \\
	 & \quad +(F_{fz}\cos\delta+F_{rz})C_r \label{F_x}  \\ 
 F_y &= F_{fy}\cos \delta + F_{fx} \sin \delta +F_{ry} + F_{fz} \sin \delta C_r \label{F_y} \\ 
 M_z &= a\left(F_{fy}\cos\delta + F_{fx}\sin\delta \right) - bF_{ry}. \label{M_z}
\end{align}
The tire force system is illustrated in Fig.\,\ref{tyre_force_diag} and is discussed in Section \ref{sub:tyre_forces}. The coefficient
$C_r$ is the tire rolling resistance coefficient and the last two terms in \eqref{F_x} and \eqref{F_y} represent the rolling resistance forces. The aerodynamic drag force $F_{ax}$ acts in  negative $x$-axis direction and is given by
\begin{equation}
F_{ax} = -\frac{1}{2} \, C_D \, \rho \, A \, u^2
\label{aero_drag}
\end{equation}
where $C_D$ is the drag coefficient.
The equations of motion (\ref{longF}), (\ref{latF}) and (\ref{yaw_bal}), expressed in terms of the elapsed arc length, are as follows
\begin{align} 
u^\prime &= S_f \dot{u} \label{u_s} \\
v^\prime &= S_f \dot{v} \label{v_s} \\
\omega_z^\prime &= S_f \dot{\bar{\omega}}_z. \label{yaw_s}
\end{align}
The angular acceleration of the Darboux (road) frame is neglected in the dynamic equations for the examples given here as its influence is negligibly small.

\subsection{Load Transfer} \label{subsub:load_transfer}
In order to compute the tire normal loads we balance the forces acting on the car normal to the road, and then balance moments around the body-fixed $y_b$ axis (see Fig.\,\ref{tyre_force_diag}). These calculations must recognise the gravitational, inertial, centripetal and aerodynamic forces acting on the car as well as the road-related gyroscopic moments. The vertical force balance gives
\begin{align}
(F_{rz}&+F_{fz})/M - n \dot{\omega}_x + g \cos \phi \cos \mu \nonumber\\
&  + (u-h \omega_y) \bar{\omega}_y -  (v+h \omega_x) \bar{\omega}_x=0 , \label{eq:load_force_balance}
\end{align}
in which the $F_{fz}$ and $F_{rz}$  are the front and rear tire normal loads, the second term 
derives from the $x$-axis component of the track's angular acceleration, 
the third term is the acceleration due to gravity, while the last two terms are centripetal accelerations.

Balancing the pitching moments around the car's mass centre gives
\begin{align}
bF_{rz}-aF_{fz} + h F_x - {I_y \dot{\bar{\omega}}_y} + (I_z-I_x) \bar{\omega}_z \bar{\omega}_x =0.  \label{eq:load_moment_balance_y}
\end{align}
The first two terms represent the pitching moments produced by the vertical tire forces, the third term is the pitching moment produced by the longitudinal force $F_x$ given in \eqref{F_x}, the fourth is the inertial moment around the car's pitch axis,
whilst the fifth term is a gyroscopic moment acting in the car's pitch direction.

\subsection{Tire Forces} \label{sub:tyre_forces}
The tire forces have normal, longitudinal and lateral components that act on the vehicle chassis at the tire ground contact points and react on the inertial frame. The rear-wheel tire force is expressed in the vehicle's body-fixed reference frame, while the front tire force is expressed in a steered reference frame; refer again to Fig.\,\ref{tyre_force_diag}. We make use of the well-known Magic Formula tire model \cite{Pacejka2006}, where these forces are a function of the normal load and the tire's longitudinal slip coefficient $\kappa$ and a lateral slip angle $\alpha$. The tire equations were also described in details in the Appendix of \cite{Limebeer20153DOptimal}. The same tire parameters are used in this work except that the peak longitudinal and lateral friction coefficients have been scaled down by 30\,\%. Following standard conventions we use
\begin{eqnarray}
\kappa &=& -\left(1+\frac{R \omega_{w}}{u_w}\right) \\
\tan{\alpha} &=& - \frac{v_w}{u_w},
\label{eq:tyreslips}
\end{eqnarray}
where $R$ is the wheel radius and $\omega_w$ the wheel's spin velocity. The quantities $u_w$ and $v_w$ are the absolute velocity components of the wheel centre in a wheel-fixed coordinate system. The front and rear tire lateral slip angles are given by
\begin{align}
\alpha_{r} &= \arctan \left( \frac{v -\dot{\psi} b}{u} \right), \\ 
\alpha_{f} &= \arctan \left( \frac{\cos\delta (\dot{\psi} a+v) -\sin\delta u}{\cos\delta u + \sin\delta (\dot{\psi} a+v)}\right).
\label{eq:4w_slip_angles}
\end{align}

\subsection{Battery Model}
\begin{figure}[ht]
\centering
\includegraphics[width=0.3\textwidth]{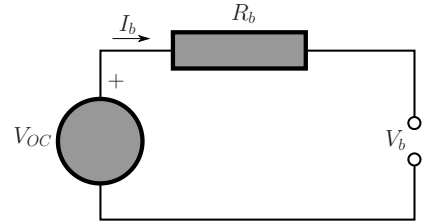}
\caption{Battery equivalent circuit model.}
\label{batteryCircuit}
\end{figure}

We will use a simple `voltage behind output resistance' model for the (lithium-ion) battery, as shown in Fig.\,\ref{batteryCircuit}.
The terminal voltage is
\begin{equation}
V_b=V_{OC}-I_b R_b,
\label{terminalVoltageEqn}
\end{equation}
where $V_{OC}$ is the battery open-circuit voltage, $I_b$ is the battery current and $R_b$ is the internal resistor.
In this convention positive $I_b$ corresponds to discharging the battery, while negative $I_b$ corresponds to charging. The unloaded battery voltage $V_{OC}$ has a dependency on the state of the charge $SoC$ as follows
\begin{equation}
V_{OC}=V^{min}_{OC} + (V^{max}_{OC}-V^{min}_{OC})SoC,
\label{ocVoltageEqn}
\end{equation}
where the $SoC$ is defined as
\begin{equation}
SoC=\frac{Q_b}{Q^{max}_b}.
\end{equation}
The maximum battery charge $Q^{max}_{b}$ is given by
\begin{equation}
Q^{max}_{b}=\frac{2 E^{max}_b}{V^{max}_{OC}+V^{min}_{OC}}
\end{equation}
where $E^{max}_b$ represents the maximum energy storage capacity of the battery. Using \eqref{terminalVoltageEqn} and \eqref{ocVoltageEqn} the power delivered, or drawn from the battery is
\begin{equation}
P_b=(V^{min}_{OC} + (V^{max}_{OC}-V^{min}_{OC})SoC- I_b R_b)I_b .
\end{equation}
Again negative $P_b$ implies charging and positive $P_b$ implies discharging.

The battery charge can be modelled by the dynamic equation
\begin{equation}
\dot{Q}_b=-I_b .
\label{batteryDynamics}
\end{equation}
The power transmission between the battery and the rear wheel requires an electric motor, which is assumed to have an efficiency factor $\mu_{em}$. The electric motor output power is thus
\begin{equation}
P_{em}=\mu_{em}^{sign(P_b)}P_b.
\label{emPowerEqn}
\end{equation}

\subsection{Engine Map}
\begin{figure}[tb]
\begin{center}
\centering
	\includegraphics[width=0.5\textwidth]{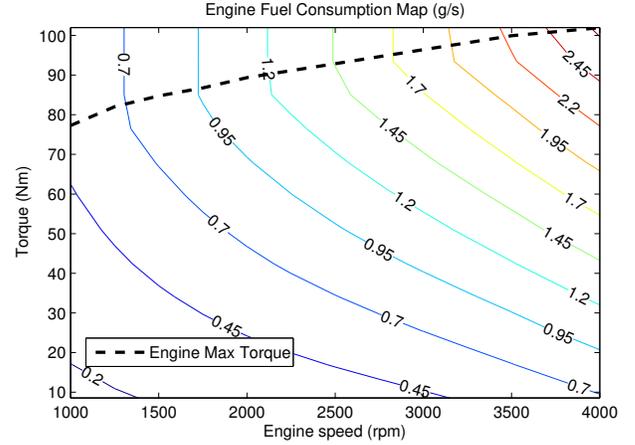}  
\end{center}
\caption{Engine fuel consumption map. The consumption rate in g/s is shown on the contours/level sets. The black dashed line is the maximum engine torque available 
as a function of speed.}
\label{engineFuelMap}
\end{figure}

The engine used in the work presented here is the 1.5\,L Prius engine with maximum power output of 43\,kW. The fuel consumption map (see Fig.\,\ref{engineFuelMap}) for this engine was obtained from the ADVISOR software \cite{ADVISOR2}. One measure of fuel efficiency is brake specific fuel consumption (BSFC), which represents the fuel mass needed to release one unit of energy. A BSFC map can be calculated from the fuel consumption map using
\begin{equation}
BSFC=\frac{\dot{m}_f(\omega_e,T_e)}{\omega_e T_e},
\label{bsfc}
\end{equation}
where $\dot{m}_f$ is the fuel-mass consumption rate and $T_e$ and $\omega_e$ are the engine torque and rotational speed respectively. We will use the BSFC to evaluate the efficiency performance of the engine. The internal combustion engine power is represented by $P_{ICE}$ and is given by
\begin{equation}
P_{ICE}=T_e \omega_e.
\end{equation}
In the optimal control calculations a quadratic multivariate polynomial was fitted to the engine map using a Linear Least Squares algorithm to speed up the fuel consumption calculation. The polynomial captures the shape of the map quite well and has an average absolute error of 2.6\% over the entire engine operating range.

\subsection{Flywheel}
A high-speed flywheel is included in the drivetrain to provide a high-power energy re-deployment capability, which complements the low-power high-energy storage capability of the battery. While batteries can store energy for a relatively long time due to their low inherent losses, flywheel storage systems suffer from high losses especially when running at high speeds.
A basic flywheel storage system can be represented by a spinning inertia with kinetic energy \begin{equation}
E_{fly}=\frac{1}{2}J_f \omega^2_f,
\end{equation}
where $J_f$ is the moment of inertia and $\omega_f$ is the flywheel's angular velocity. The dominant losses in the flywheel come from the friction in the bearings. These losses are modelled using the empirical relationship given on  pg.\,147 of \cite{Doucette2013}
\begin{equation}
P_{loss}=2 \times 10^{-7} \omega^2_f + 0.0151 \omega_f + 4.0577, \label{FlywheelLosses}
\end{equation}
where $\omega_f$ is given in $rpm$. The flywheel dynamics are described by
\begin{equation}
\dot{E}_{fly}=-P_{fly}-P_{loss} .
\label{flywheelDynamics}
\end{equation}
There are also losses in the continuous variable transmission (CVT) that can be lumped into an efficiency factor $\mu_{CVT}$.

\subsection{Power Transmission}
The power transmitted to the rear wheel can be modelled by the constraint
\begin{equation}
P_{ICE}  + P_{em} + \mu_{CVT}^{sign(P_{fly})} P_{fly}-(F_{rr x}+F_{rl x})u \geq 0,
\label{powerConstraint}
\end{equation}
which ensures that the power delivered to the back wheels never exceeds the combined power delivery capability of the internal combustion engine, the battery and the flywheel. If the rear-wheel tire force is positive, the vehicle is being driven, and the sum of powers delivered by the engine, the flywheel CVT and electric motor will match the mechanical power delivered to the rear wheel. Under braking, rear-wheel tire force is negative, and the mechanical power at the back wheels is used to charge the battery and/or the flywheel, or else it is dissipated as heat.

At $P=0$, $\mu_{CVT}^{sign(P_{fly})}$ is discontinuous and hence must be approximated using a smooth function. The approximation used in this study is
\begin{align}
\mu^{sign(P)} &\approx 0.5 \mu\left(1+\tanh(\varrho P)\right) + \frac{0.5}{\mu}(1+\tanh(-\varrho P)),
\label{efficiencyApprox}
\end{align}
in which $\varrho$ is a constant.  As $\varrho$ is increased, the approximation (\ref{efficiencyApprox}) approaches $\mu^{sign(P)}$.

\section{Optimal Control} \label{OptimalControl} The minimum fuel problem can be formulated as an optimal control problem that is now described.

The system is described by a set of equations of the form 
\begin{equation}
x^\prime (s) =\bm{f}(\bm{x}(s),\bm{u}(s)),
\end{equation}
in which the state-vector is given by
\begin{equation}
\bm{x}=
\begin{bmatrix}
n &\xi & u& v& \omega & E_{fly} & Q_b
\end{bmatrix}^T.
\end{equation}
The associated differential equations are given by \eqref{n_rate_s}, \eqref{xi_rate_s}, \eqref{longF}, \eqref{latF}, \eqref{yaw_bal}, \eqref{batteryDynamics} and \eqref{flywheelDynamics}, respectively.
The control vector is given by
\begin{equation}
\bm{u} = 
[\begin{matrix} \delta & k_{f} & k_{r} & F_{fz} &  F_{rz} & \omega_e & T_e & P_{fly} &  I_b\end{matrix}]^T.
\end{equation}
The problem is also subject to the path constraints \eqref{eq:load_force_balance}, \eqref{eq:load_moment_balance_y}, \eqref{powerConstraint} and bounds on the states and controls. There is also a constraint on maximum engine
torque available as shown in Fig.\,\ref{engineFuelMap}.

The minimum-fuel performance index to be minimised is given by 
\begin{equation}
J =  \int\limits_{s_0}^{s_f}  S_f(s)\dot{m}_f \, ds.
\label{eqnLagraneCost}
\end{equation}
The arrival time constraint can be written as an integral constraint as below
\begin{equation}
\int\limits_{s_0}^{s_f} S_f(s) ds \, \leq  \,T,
\end{equation}
where $T$ is the arrival time.

In practice, however, to avoid jerky controls and singular arcs, we actually control $\dot{\bm{u}}$ and minimise
\begin{equation}
J_{mod} =  \int\limits_{s_0}^{s_f}  S_f(s)(\dot{m}_f +\dot{\bm{u}}^T \bm{R} \dot{\bm{u}}) \, ds
\end{equation}
in which $\bm{R}$ is an appropriate weighting matrix. We also impose slew-rate limits on controls by placing hard constraints on $|\dot{\bm{u}}|$.

\subsection{Numerical Optimal Control} \label{optimalControl}
An optimal control problem formulation general enough for our purposes is of Lagrange form. The aim is to determine states $\bm{x}(\tau) \in \mathbb{R}^n$, controls $\bm{u}(\tau) \in \mathbb{R}^m$ and static parameters  $\bm{p}(\tau) \in \mathbb{R}^q$ which minimise a cost functional
\begin{equation}
J = \frac{t_f-t_0}{2} \int\limits_{-1}^{+1}g[\bm{x}(\tau),\bm{u}(\tau),\tau, t_0,t_f,\bm{p}] d\tau
\label{eqnCostFucntional}
\end{equation}
subject to the state dynamics,
\begin{equation}
 \frac{d\bm{x}}{d\tau} =\frac{t_f-t_0}{2} \bm{f}[\bm{x}(\tau),\bm{u}(\tau) ,\tau,t_0,t_f,\bm{p}],  
\label{eqnStateDynamics}
\end{equation}
path constraints
\begin{equation}
 \bm{c}_{min} \leq \bm{c}[\bm{x}(\tau),\bm{u}(\tau), \tau,t_0,t_f,\bm{p}] \leq  \bm{c}_{max} \in \mathbb{R}^r,
\end{equation}
and boundary conditions
\begin{equation}
 \bm{b}_{min} \leq \bm{b}[\bm{x}(-1),\bm{x}(+1), t_0,t_f,\bm{p}] \leq  \bm{b}_{max} \in \mathbb{R}^s.
\end{equation}
The normalised optimisation interval $\tau \in [-1, 1]$ can be transformed into the general interval $t \in [t_0 ,t_f ]$ using the affine transformation $t = (t_f - t_0)\tau /2 + (t_f + t_0)/2$. 

The pseudo-spectral numerical optimal control solver GPOPS-II \cite{Patterson2013}, which is based on the Legendre-Gauss-Radau (LGR) collocation scheme, was used to solve the minimum lap time problem in this paper. In this scheme the state is approximated using a Lagrange polynomial of order $N$
\begin{equation}
 \bm{x}(\tau)\thickapprox \bm{X}(\tau) = \sum_{i=1}^{N+1} \bm{X}_i L_i(\tau) 
\end{equation}
where
\begin{equation}
L_i (\tau) = \prod_{\substack{
             j = 1\\
             j \ne i }}^{N+1}
       \frac{\tau - \tau_j}{\tau_i - \tau_j} , i=1, \ldots, N+1.
\end{equation}

The state-derivative approximation is thus given by
\begin{equation}
 \dot{\bm{x}}(\tau)\thickapprox \dot{\bm{X}}(\tau) = \sum_{i=1}^{N+1} \bm{X}_i \dot{L_i}(\tau). 
\end{equation}

Collocating the state dynamics at $N$ LGR points gives
\begin{align}
 \dot{\bm{X}}(\tau_j)\ = \sum_{i=1}^{N+1} \bm{X}_i \dot{L_i}(\tau_j) = \sum_{i=1}^{N+1}  \bm{X}_i D_{ji},~j=1, \ldots,N,
 \label{eqnStateDeriv}
\end{align}
where $D_{ji}= \dot{L_i}(\tau_j)$ are the $N \times (N+1)$ elements of the LGR differentiation matrix. Note that $\tau_{N+1}$ is a non-collocated point. Using \eqref{eqnStateDeriv} we can discretise \eqref{eqnStateDynamics} and essentially transform the state dynamics given by ordinary differential equations into algebraic constraints.

The optimal control problem can then be approximated by an NLP problem. The cost function of this NLP is obtained by approximating the cost functional \eqref{eqnCostFucntional} using LGR quadrature. The NLP problem can be described by the task of finding $\bm{X}_i$'s, $\bm{U}_i$'s and $\bm{p}$ which minimise
\begin{equation}
J \thickapprox \frac{t_f-t_0}{2} \sum_{i=1}^{N}w_i g[\bm{X}_i,\bm{U}_i,\tau_i, t_0,t_f,\bm{p}] d\tau
\end{equation}
in which the $w_i$'s are the quadrature weights \cite{Abramowitz1965},
subject to the following constraints
\begin{align}
  \sum_{i=1}^{N+1}  \bm{X}_i D_{ji} &=  \frac{t_f-t_0}{2} \bm{f}[\bm{X}_i,\bm{U}_i,\tau_i, t_0,t_f,\bm{p}], ~j=1, \ldots,N;
\end{align}

\begin{align}
 \bm{c}_{min} \leq \bm{c}[\bm{X}_i,\bm{U}_i, \tau_i,t_0,t_f,\bm{p}] \leq  \bm{c}_{max} \newline ,~i=1, \ldots,N;
\end{align}

\begin{align}
 \bm{b}_{min} \leq \bm{b}[\bm{X}_1,\bm{X}_N, t_0,t_f,\bm{p}] \leq  \bm{b}_{max}.
\end{align}

For clarity of exposition, the description provided here is for a single-interval pseudo-spectral method (global collocation). GPOPS-II uses a mesh and polynomial degree refinement scheme \cite{Darby2011} so that error reduction can be achieved in the presence of non-smooth problem features. The extension to multiple segments is straightforward, with the only requirement being the need to enforce continuity between each mesh interval.

The transcribed NLP problem is typically large but sparse. The IPOPT\cite{Biegler2008} software library (based on interior point methods) was used to solve the NLP problem. Automatic Differentiation was used to provide IPOPT with accurate and computationally efficient first and second order derivatives \cite{Weinstein2014}. 

\section{Results} \label{Results}
\begin{figure}[tb]
	\centering
		\includegraphics[width=0.5\textwidth]{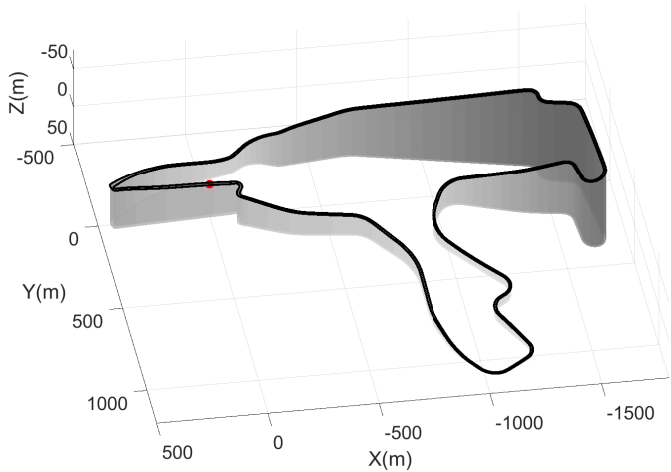}
		\caption{Circuit de Spa-Francorchamps used in the presented study. Start\slash Finish line (SF) is marked by the blue dot.}			
		\label{Spa3d}
\end{figure}

\begin{figure}[tb]
	\centering
		\includegraphics[width=0.45\textwidth]{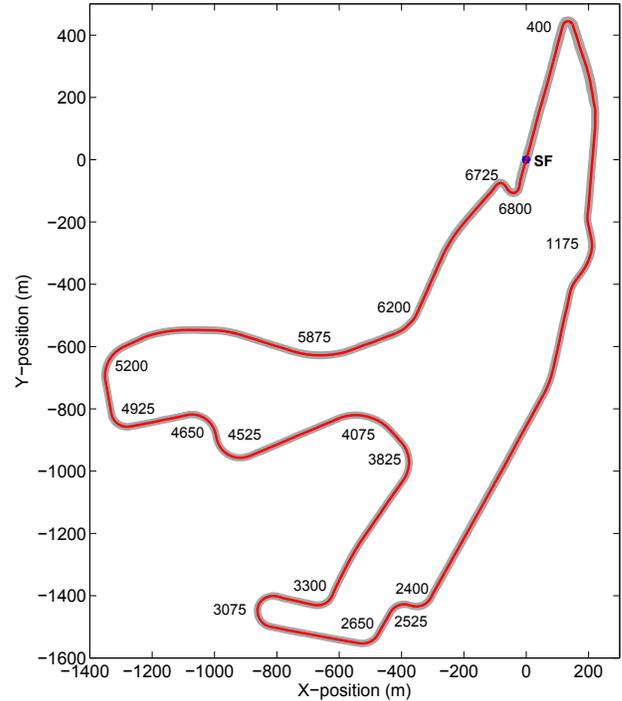}
		\caption{The vehicle's path with the distances to the various corners from the SF line (in meters).}			
		\label{RacingLine}
\end{figure}

In order to illustrate the concepts described, the motor-racing Circuit de Spa-Francorchamps will be used as an exemplar track. This track is approximately 7\,km long with an elevation change of approximately 110\,m.
The path is restricted to the neighbourhood of the track centre line in order to make comparisons easier; this was achieved by constraining 
the state $n$, which is the perpendicular distance from the car mass centre to the track centre line, to be `small'. The three-dimensional track, as well as a two-dimensional projection on to a ground plane are shown in Fig.\,\ref{Spa3d}.
The route that the car is constrained to take together with the corner distances are shown in Fig.\,\ref{RacingLine}. 
This figure will be useful in analysing the results to be presented later. In all the simulations described here, the vehicle will start at rest from the start-finish line (SF) (in practice the vehicle will start from a low speed in order to keep its reciprocal well defined), and will complete the circuit such that the combustible fuel usage is minimised. The car will come to a standstill at the end of its journey. The vehicle parameters are summarised in Table\,\ref{vehicleParameters} given in the Appendix.

\begin{figure}[tb]
	\centering
		\includegraphics[width=0.5\textwidth]{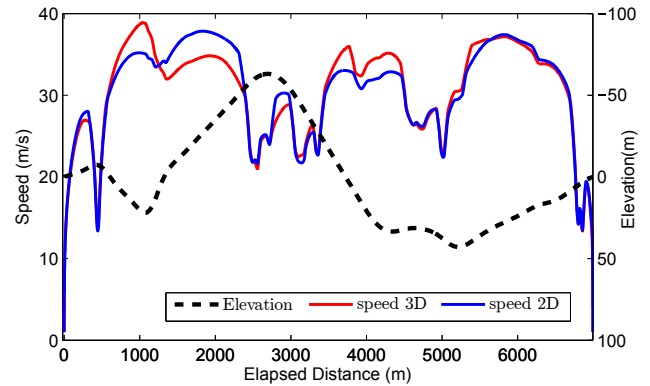}
		\caption{Speed profile for the three dimensional (red solid line) and two dimensional (blue solid line) track for the case when arrival time was set to 240\,s. The black dashed line shows the road elevation heights. }			
		\label{speedProfile}
\end{figure}

\begin{figure}[tb]
	\centering
		\includegraphics[width=0.5\textwidth]{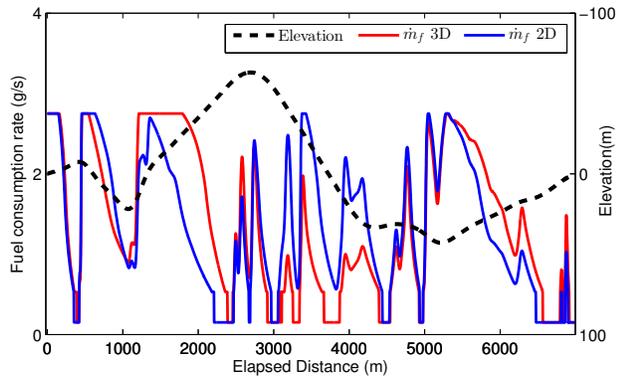}
		\caption{Engine fuel consumption rate in g/s for the three dimensional (red solid line) and two dimensional (blue solid line) track for the case when arrival time was set to 240\,s. The black dashed line shows the road elevation profile. }			
		\label{fuelRate}
\end{figure}

The minimum fuel problem was solved for a journey time of 240\,s on the two-dimensional and three-dimensional track descriptions. These comparitive calculations are used to quantify the effects of three dimensionality.
The optimal speed profile of the full hybrid vehicle for both the two and three dimensional tracks, along with the track elevation changes, are depicted in Fig.\,\ref{speedProfile}.
One can see that the track is on a slight incline when the car starts its journey.
This results in the car's speed at 200\,m from the SF line being 
slightly higher on the 2D track. After the hairpin bend at 400\,m 
the track falls away and the predicted speed on the 3D track exceeds that of the 2D track model.
The 3D track speed advantage is then `given back' as the car enters the uphill section between 1100\,m and 2400\,m. At the start of the incline at 1100\,m, one also sees an increase in the fuel consumption rate as shown in Fig.\,\ref{fuelRate}.
\begin{figure*}[tb]
	\centering
		\includegraphics[width=1.0\textwidth]{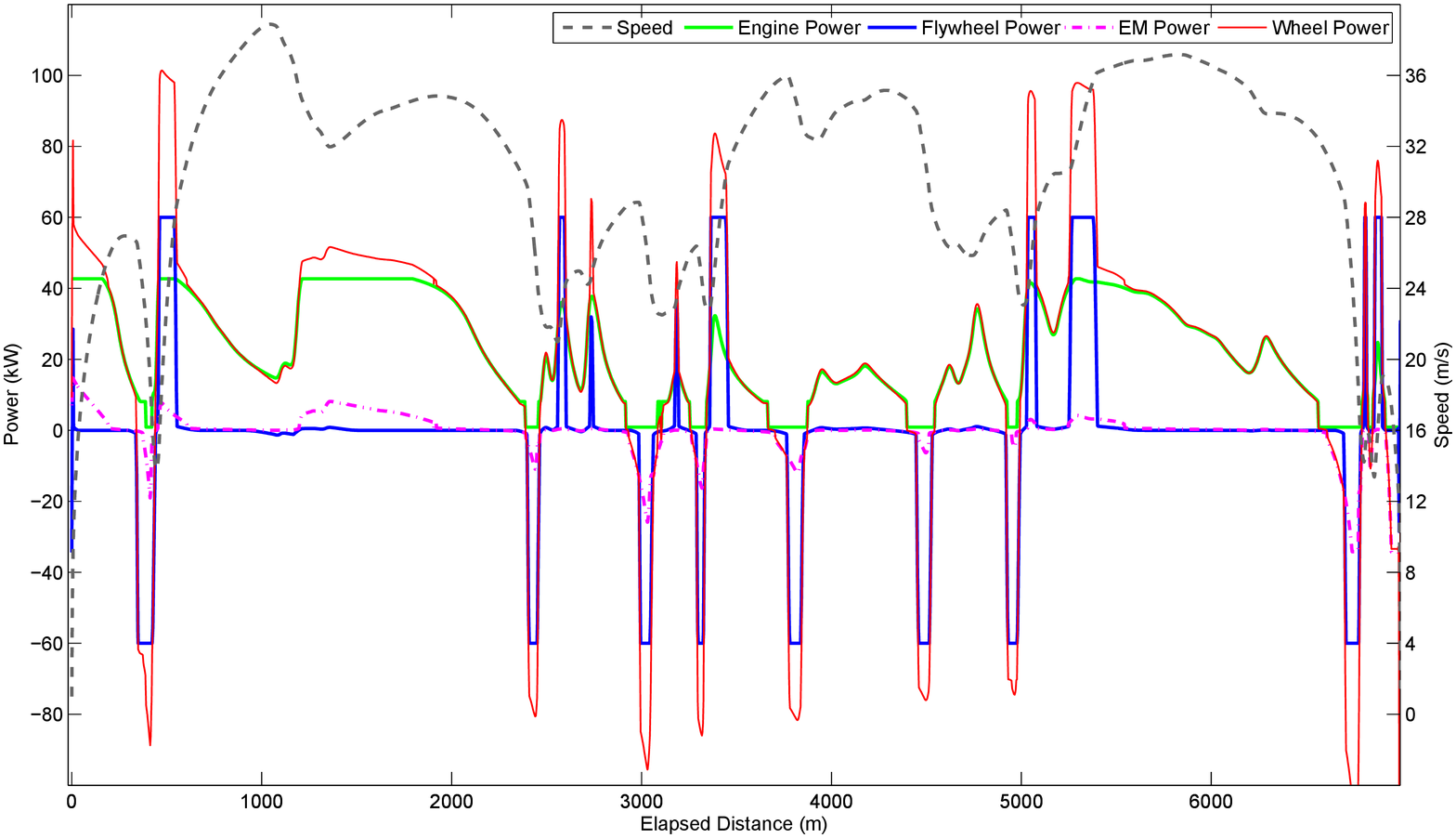}
		\caption{Power strategy for the vehicle with engine, flywheel and battery when the arrival time was set to 240\,s. The power produced by the engine is shown in green, flywheel power in blue, electric motor power in magenta and the rear wheel mechanical power is shown in red. The grey dashed line represents the speed profile. }			
		\label{powerPlots}
\end{figure*}

The vehicle's fuel consumption rate is high initially in order for it to accelerate from rest. As the vehicle approaches the hairpin bend at 400\,m, the throttle is eased off at approximately 150\,m, with the brake applied at approximately 200\,m (see also Fig.\,\ref{powerPlots}). The full-lap fuel usage for the 2D case was 321.2\,g, while that in the 3D case was 320.5\,g. Although this difference is small, Fig.\,\ref{fuelRate} shows that the fuel consumption strategy is different in the two cases; more fuel is used ascending hills, while fuel is saved coming down them.

Fig.\,\ref{powerPlots} shows the power management strategy of the vehicle. An inspection of this plot brings a number of issues to light: (i) the battery power delivered and absorbed is relatively small compared to the flywheel; (ii) the controller tries to utilise the stored flywheel energy as fast as possible in order to avoid energy dissipation resulting from the rotational losses; (iii) the battery stored energy is only used when that in the flywheel has been depleted; (iv) often times the power delivery from the flywheel exceeds that of the engine; (v) the battery is discharged at relatively few isolated points on the circuit, and only when the flywheel energy has been depleted.

\begin{figure}[tb]
	\centering
		\includegraphics[width=0.49\textwidth]{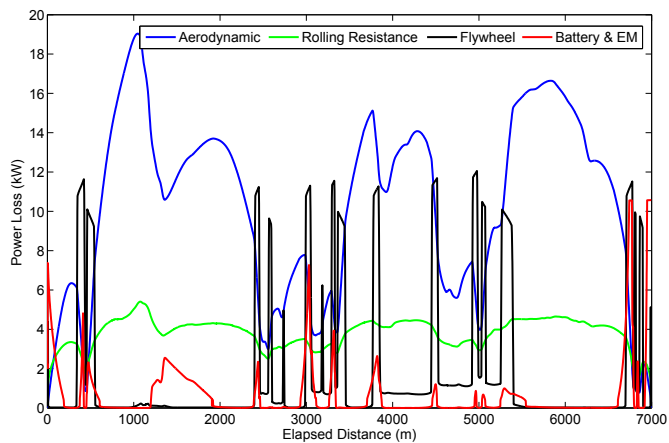}
		\caption{Power losses as a function of elapsed distance from the start-finish line. The aerodynamic losses are shown in blue, the rolling resistance losses are shown in green, the flywheel losses are shown in black, and the battery and motor losses in red.}			
		\label{LossesDistance}
\end{figure}

The power losses associated with various dissipation mechanisms are shown in Fig.\,\ref{LossesDistance}. Predominant are the aerodynamic drag power losses that are proportional to the cube of the forward speed; see \eqref{aero_drag}. The rolling resistance loss is given by the product of forward speed and the rolling resistance force, which is also relatively large in this case. The flywheel losses are the sum of internal rotational losses (\ref{FlywheelLosses}) and the CVT losses. The peaks which appear in flywheel loss plot are substantially due to transmission losses. The positive gradients appearing on some of the peaks are associated with the increasing losses as the flywheel is charged. The negative gradients, on the other hand, are associated with decreasing losses as the flywheel is discharged. The low loss regions, such as the one which appears between 3900\,m and 4300\,m, are due to self-discharging. The battery and electric motor losses come from the battery's output resistance and motor losses.

\begin{figure}[tb]
	\centering
		\includegraphics[width=0.5\textwidth]{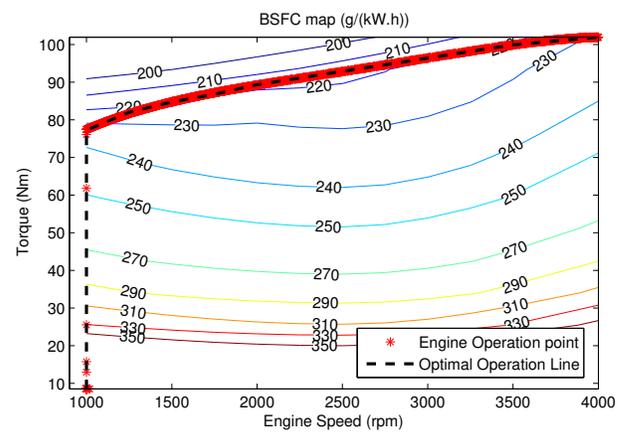}
		\caption{Engine operation points on the BSFC map. The optimum operation line is shown in dashed black line. }			
		\label{BSFC}
\end{figure}

The brake specific fuel consumption map \eqref{bsfc} is shown in Fig.\,\ref{BSFC}. 
It can be seen that the optimal controller tries to operate the engine as close as possible to the optimum efficiency line; the controller chooses the engine speed and torque so that maximum mechanical energy is extracted per unit of fuel mass burnt.

\begin{figure}[tb]
	\centering
		\includegraphics[width=0.5\textwidth]{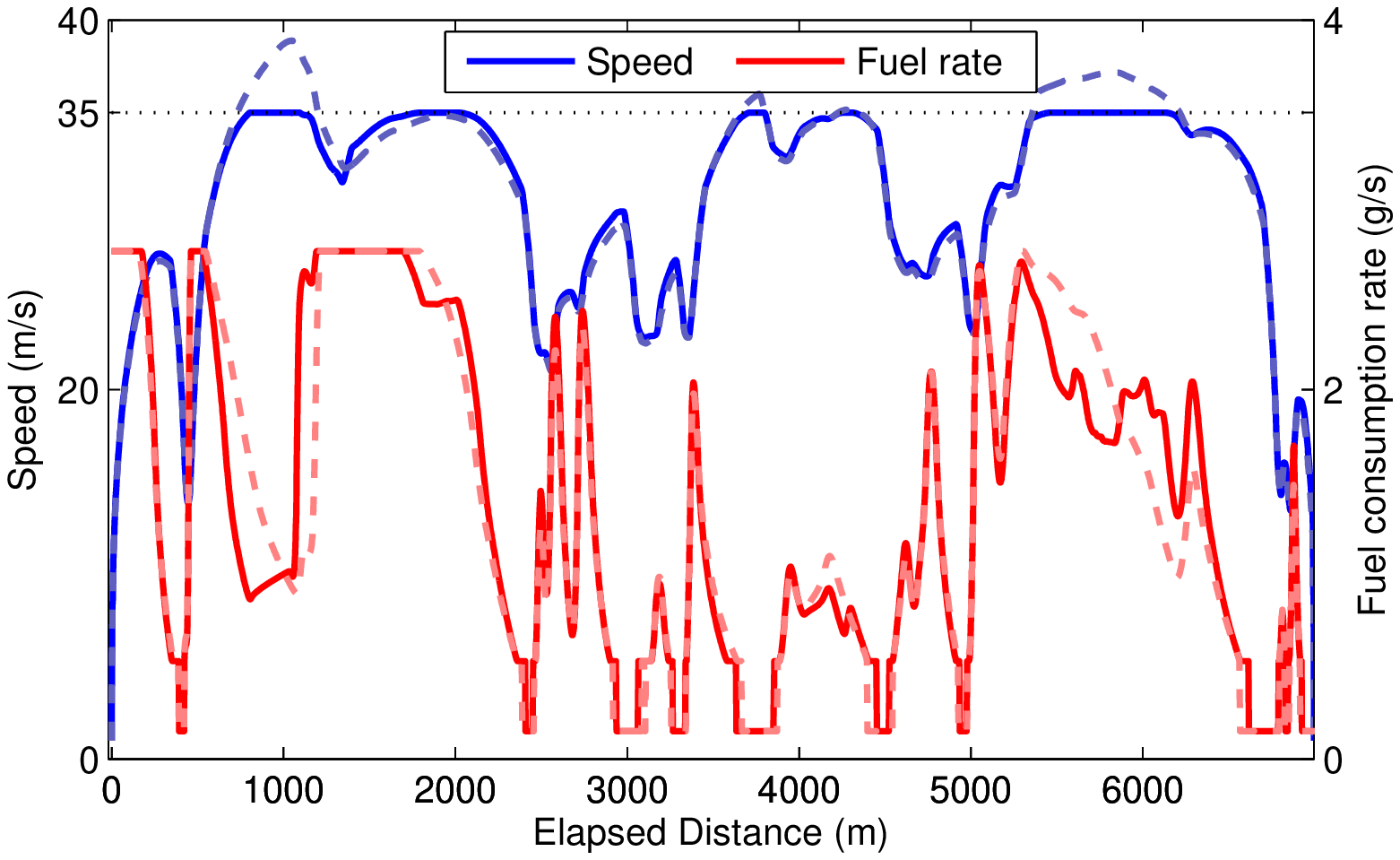}
		\caption{Speed (solid blue) and fuel consumption rate (solid red) when a speed-limit of 35\,m/s is imposed. The dashed lines correspond to the no speed-limit case. }			
		\label{FigSpeedLimitCase}
\end{figure}

The results presented thus far can be expanded to include real-world influences such as speed limits, enforced stops, changes in the road-surface conditions and wind gusting. One approach to the solution of these problems is to set them up as multi-phase optimal control problems \cite{Rao2010}. In this formalism each phase can contain new models and/or new model parameters, new inputs, and new constraints. New parameters might include down-graded tyre parameters that particularise degraded road surface conditions, new inputs might include wind-related disturbances, and new constraints might include such things as speed limits and traffic controls. Fig.\,\ref{FigSpeedLimitCase} demonstrates how the speed and fuel consumption rate vary when a 35\,m/s speed-limit is imposed. At an elapsed distance of 1000\,m it is evident that the fuel consumption rate reduces (below the unrestricted speed case) on entry to the speed-restricted section of road, and then increases above the unrestricted speed case in order to make up for the time lost. Similar variations in fuel consumption can be observed on the 5300\,m to 6200\,m  road section.

\begin{figure}[tb]
	\centering
		\includegraphics[width=0.5\textwidth]{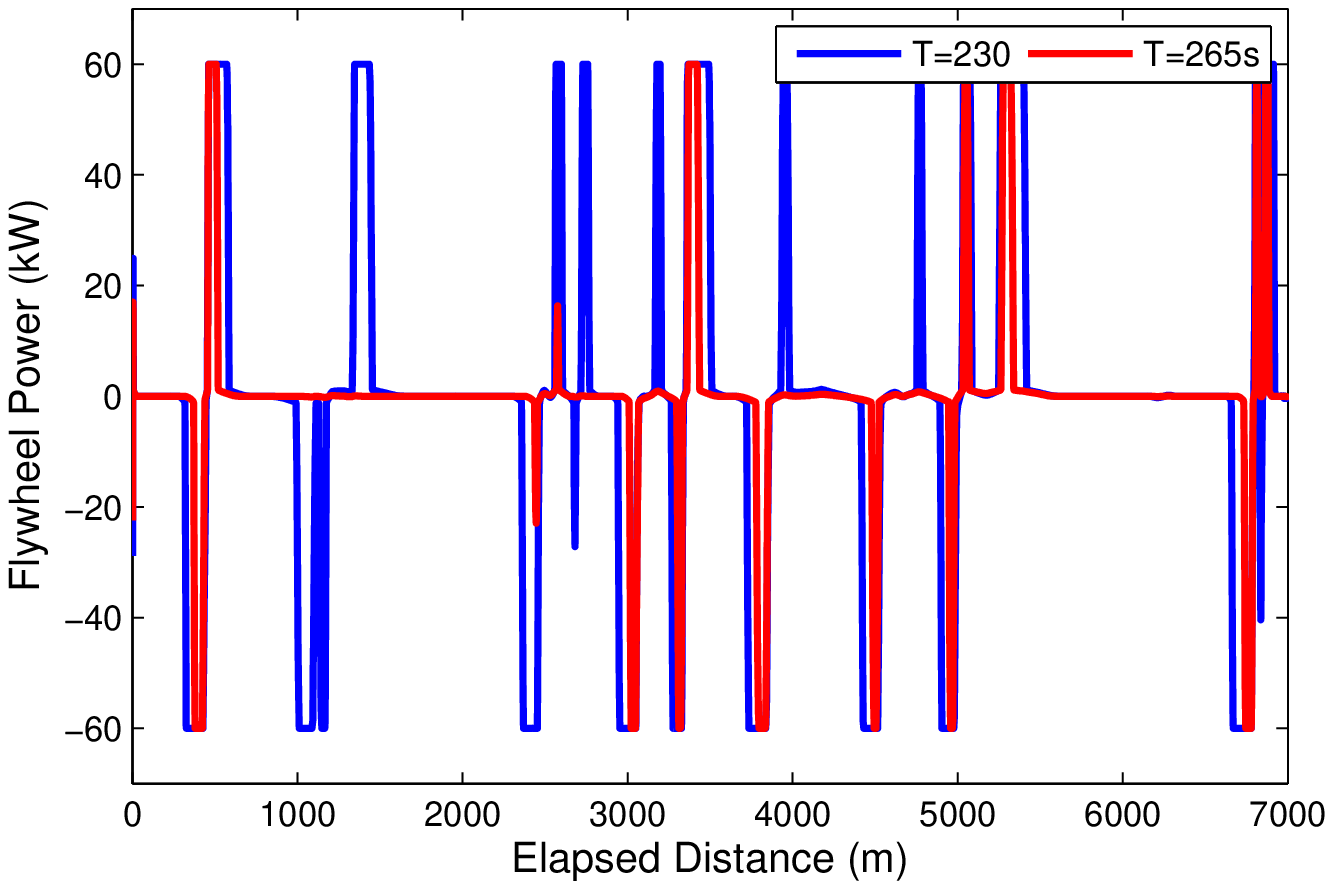}
		\caption{Flywheel power for different journey times for the vehicle with no battery. }			
		\label{flywheelPower}
\end{figure}

\begin{figure}[tb]
	\centering
		\includegraphics[width=0.5\textwidth]{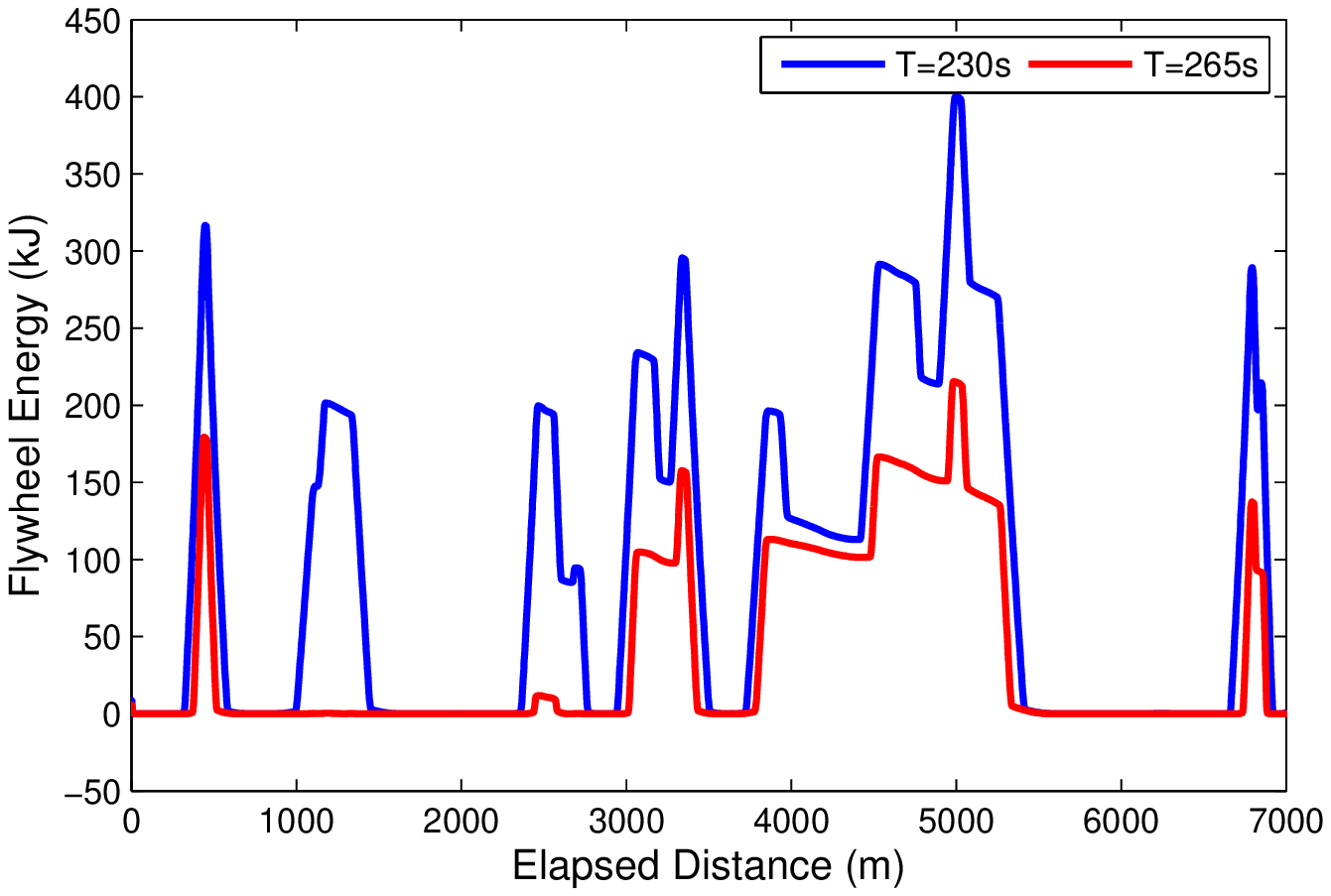}
		\caption{Flywheel energy for different journey times for the vehicle with no battery. }			
		\label{flywheelEnergygLine}
\end{figure}

It will be shown that the flywheel can provide significant fuel savings, especially when the journey times are low, and when the vehicle is required to complete the route aggressively. 
In these cases there will be many braking regions that will regenerate energy back into the flywheel that can then be quickly redeployed to accelerate the vehicle. For longer journey times the vehicle can complete the circuit at low speed with little or no braking. These ideas are illustrated in Figs.\,\ref{flywheelPower} and \ref{flywheelEnergygLine}, which show the power and energy stored in the flywheel, respectively, for a vehicle with engine and flywheel only.
For a journey time of 230\,s, the flywheel frequently reaches high levels of stored energy. In contrast, when the journey time is increased to 265\,s, the braking regions become less frequent, and lighter, resulting in fewer opportunities to scavenge energy to recharge the flywheel.
Flywheel self-discharging power losses are evident as negative gradients on the flywheel energy peaks in Fig.\,\ref{flywheelEnergygLine}. 

\begin{figure}[tb]
 	\centering
 		\includegraphics[width=0.49\textwidth]{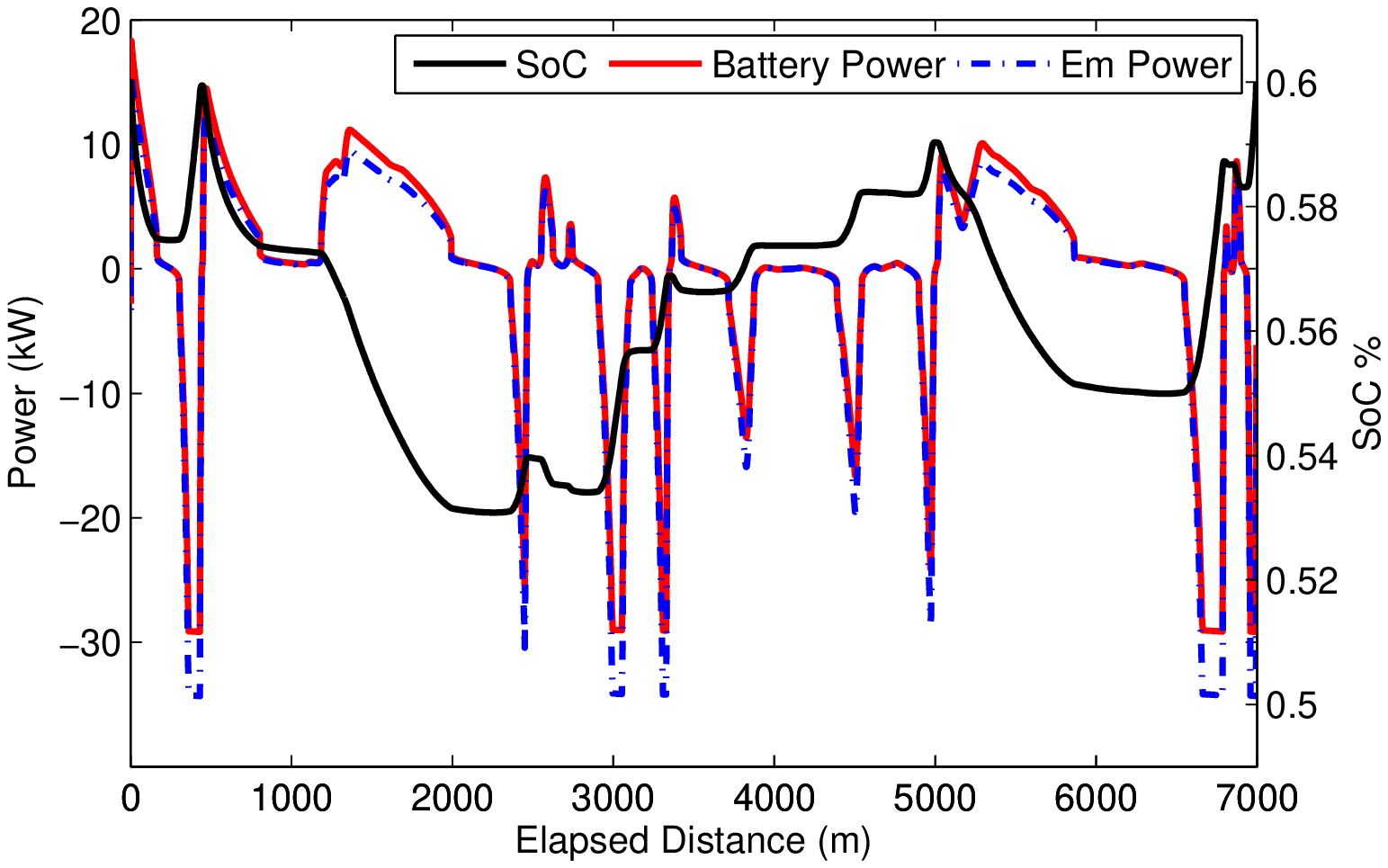}
 		\caption{The battery power and state of charge along with the power at the electric motor shaft. The journey time was 250\,s. }			
 		\label{batterySoc}
 \end{figure}

In order to analyse the battery management strategy, the vehicle with ICE and battery is considered in  Fig.\,\ref{batterySoc} for a journey time of 250\,s. 
The battery state of charge at the start and end of the lap is constrained to be 60\,\%. 
The battery is discharged in the acceleration regions for power and the re-charged in braking phases.
It is evident that unlike the flywheel, the discharge does not happen in an `on-off' fashion. This is because the power losses in the battery are proportional to square of the current.
This results in the optimal controller spreading the power delivery over longer time periods in order to maintain high efficiency. In the short charging phases all the available braking power is absorbed. In the regeneration phases, the electric motor losses result in less power than the power available at the wheels being delivered to the battery.  
These efficiency losses are also evident in the power assist mode when part of battery power is lost on its way to the wheels, and explains why the electric motor power is higher in braking and lower in acceleration.

\begin{figure}[tb]
	\centering
		\includegraphics[width=0.45\textwidth]{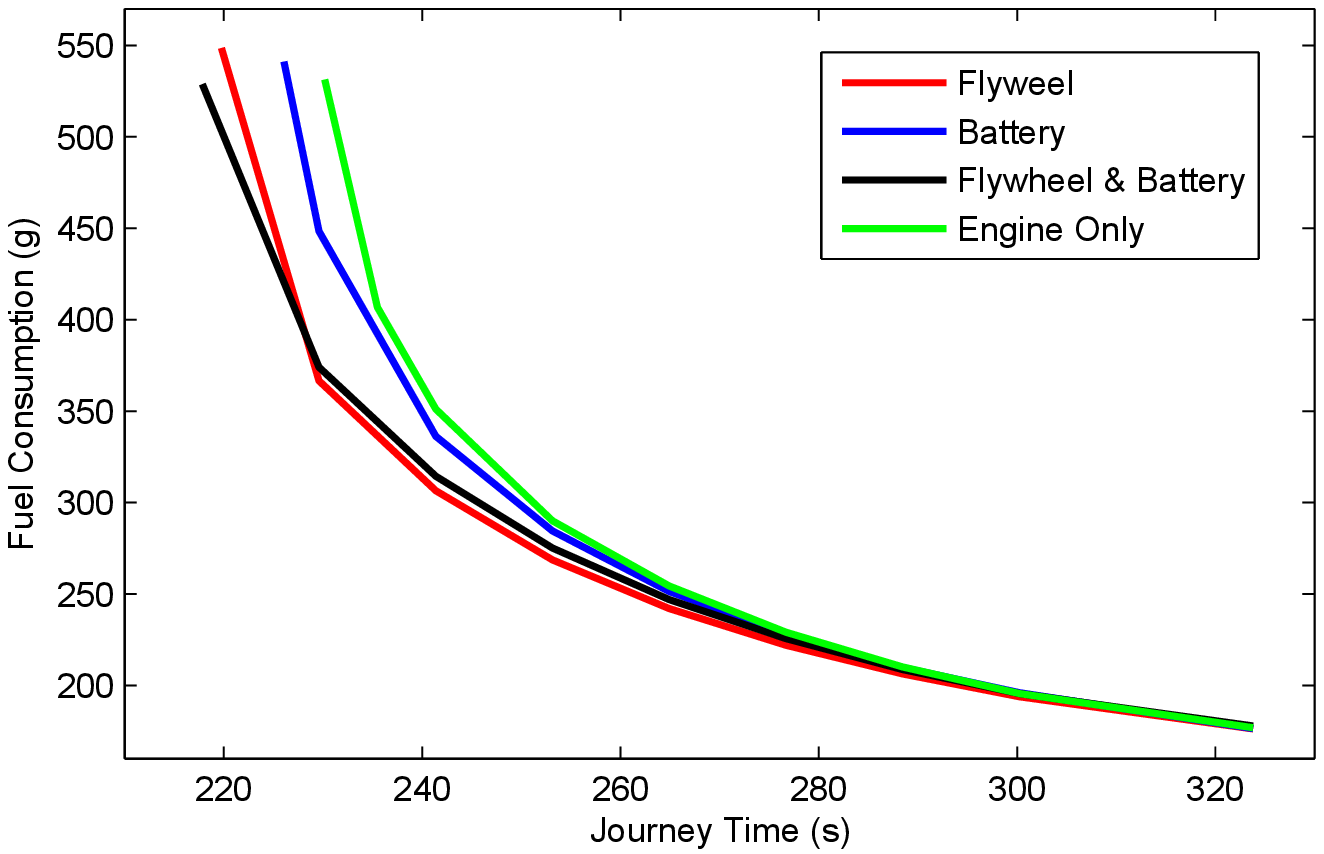}
		\caption{Engine fuel usage against the time-of-arrival constraint for completing the circuit for different storage configurations. The case with engine and flywheel is shown in red, engine and battery in blue, combination of three in black and engine only in green.}			
		\label{fuelvsJourneyTime}
\end{figure}

In a final study, the minimum fuel problem was solved for a number of arrival times for four different energy storage combinations; the results are summarised in Fig.\,\ref{fuelvsJourneyTime}. 
As one would expect, at lower journey times, when vehicle has to be driven more aggressively, the fuel consumption increases. Furthermore, at higher speeds the aerodynamic drag losses increase.
 The vehicle with both the flywheel and battery offers maximum advantage at low journey times.
However, as the journey time increases, the vehicle with engine and flywheel performs the best as it does not have to carry the additional 100\,kg of battery load. With increasing journey times the fuel usage drops significantly and the benefit of using an axillary storage system also diminishes, as energy regeneration opportunities decreases. The minimum journey times possible for the vehicle with engine only, engine and battery only, engine and flywheel only, and engine, battery and flywheel were 230.09\,s, 225.76\,s, 219.67\,s and 217.35\,s respectively. This highlights the power boost capability of the flywheel in aggressive manoeuvrings. 
 
\subsection*{Some computational details} In this study the optimal control solver was initialised with 100 mesh segments. The number of collocation points was allowed to vary between 4 to 10  per mesh segment. The error tolerance across all meshes was $10^{-3}$ and IPOPT tolerance was set to $10^{-7}$. For the vehicle with both battery and flywheel, when arrival time was set to 240\,s, the error tolerance was reached after 18 mesh adaptation iterations. The final mesh had a total of 367 mesh segments and 2861 collocation points. The whole problem took under 2 hours to solve on an 8 core 3.5 GHz computer. For the engine-only case the error tolerance was reached with a similar number of collocations and meshes. However, the total solution time was only 20 minutes as only 10 mesh adaptation iterations were required and each iteration was faster to complete. The solution time is sensitive to the problem set-up. For example, increasing the size of $\varrho$ in \eqref{efficiencyApprox} slows down convergence as the problem becomes `less continuous'. Avoiding the use of look up tables in the cost function speeds up the algorithm significantly. In general, even though the direct pseudo-spectral method shows great robustness, careful attention is required in problem set-up when fast solution speeds are desired.

\subsection*{Global Optimality}
As with any gradient based optimisation method, the optimal solution obtained from the algorithm might be a local rather than the global minimum. It is therefore important to quantify the sensitivity of the optimal solution to the initial guess. In this study, the problem was initialised using a set of sensible constant values across the entire solution space. A series of cases were then considered to examine whether the solution obtained can be trusted to be a global optimum. In all the cases, the vehicle was equipped with a battery and flywheel and the arrival time was set to 240\,s.
 
Firstly, the initial state and control guesses were chosen constant, and given in terms of the bounds by $\bm{x}_{0}=\bm{x}_{min}+(\bm{x}_{max}-\bm{x}_{min})w$ and $\bm{u}_{0}=\bm{u}_{min}+(\bm{u}_{max}-\bm{u}_{min})w$, with $w$ a weighting between zero and one.
The optimal fuel usage for these cases are shown in Table \ref{initialGuessTable} as `Linear x\%'. For the 10\%, 20\% and 80\% cases, the solution did not converge (marked as DNC). However, for all the other cases the solutions were essentially identical with the differences being less than the mesh tolerance error of $10^{-3}$.
\begin{table}[tb]
\caption{Optimal fuel usage for different initial guesses. DNC stands for `Did not converge'. }
\label{initialGuessTable}
\begin{center}
\begin{tabular}{|c|c|c|c|c|c|c|}
\hline
Initial Guess & Fuel (g) & Initial Guess & Fuel (g)  \\ 
\hline
Linear  10\%	&DNC	& Linear  20\% 	 &  DNC  	   \\ \hline
Linear  30\%	&320.53	& Linear  40\% 	 &  320.52	   \\ \hline
Linear  50\%	&320.53 	& Linear  60\%	 &  320.59     \\ \hline 
Linear  70\%	&320.52 	& Linear  80\%	 &  DNC     \\ \hline 
Random 1	& 320.51 	& Random 2	 &  320.59   \\ \hline
Random 3	& 320.59 	& Random 4 	 &  320.60      \\ \hline
Random 5	& 320.59 	& Random 5 	 &  320.59    \\ \hline 
Random 6	& 320.59 	& Random 7 	 &  320.59    \\ \hline 
Random 8	& 320.58 	& Random 9 	 &  320.57   \\ \hline 
Battery min	& 320.56 	& Battery max 	 &  320.53    \\ \hline 
\end{tabular}
\end{center}
\end{table}

\begin{figure}[tb]
	\centering
		\includegraphics[width=0.5\textwidth]{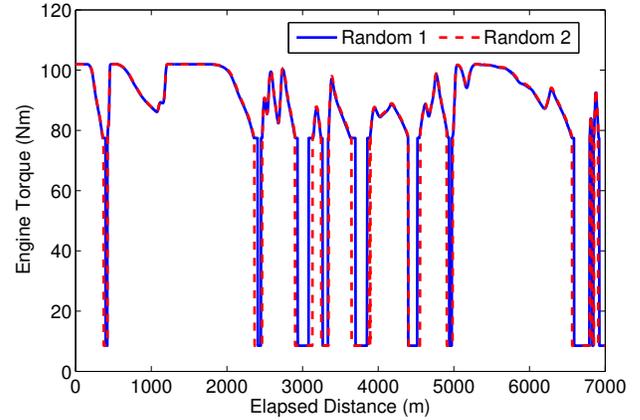}
		\caption{Optimal engine torques for the Random 1(sold blue) and Random 2(dashed red) initial guesses.  }			
		\label{intialGuessesEngineTorque}
\end{figure}

In an alternative test, a random value of between 20\% and 80\% of the total bound range, for each state and control, was chosen as the initial guess. 
This experiment was repeated 10 times and the results are shown as `Random n' in Table \ref{initialGuessTable}. All the runs converged successfully to the same solution. Finally, the problem was initialised with the lowest charging current and $SOC$, and then with highest discharge current and $SOC$ for the battery. The results are shown as `Battery min' and `Battery max'. Identical solutions were obtained once more, demonstrating that the optimal control problem has a good `radius of convergence' when solved using a direct pseudo-spectral method.
To ensure that the same minimal cost was not obtained from different trajectories the cases with the highest cost difference were selected (Case Random 1 and Random 2). The state and controls for the two cases where then compared. The worst case difference was found to be in the Engine Torque, $T_e$, as plotted in Fig.\,\ref{intialGuessesEngineTorque}. As can be seen, the trajectories are nearly identical, lending confidence to the idea that the solution obtained is indeed globally optimal.

\section{Conclusions} \label{Conclusion}
We have presented a novel method of evaluating the fuel-consumption performance of hybrid vehicles
with multiple secondary energy sources. Rather than utilising standardised driving cycles, we use a specific route and then drive the vehicle so as to we reach the destination within a given journey time, while controlling the vehicle in order to minimise the combustible fuel consumption.
In this framework driving aggressiveness can be systematically controlled by changing the arrival time. Another thrust of this work was to quantify the effectiveness of flywheels and batteries in reducing fuel consumption. Finally, the minimum fuel usage strategy for the hybrid vehicle over a three dimensional route was evaluated and the effects of different combinations of auxiliary energy storage systems were studied.

In future work it might be interesting to consider the effect of route optimisation. Adding traffic information to the model will also make the simulations more realistic. Component sizing of the combustion engine and auxiliary storage system can be formulated easily in the optimal control problem framework by including static parameters.
\bibliographystyle{IEEEtran}
\bibliography{}

\begin{thebibliography}{10}
\providecommand{\url}[1]{#1}
\csname url@samestyle\endcsname
\providecommand{\newblock}{\relax}
\providecommand{\bibinfo}[2]{#2}
\providecommand{\BIBentrySTDinterwordspacing}{\spaceskip=0pt\relax}
\providecommand{\BIBentryALTinterwordstretchfactor}{4}
\providecommand{\BIBentryALTinterwordspacing}{\spaceskip=\fontdimen2\font plus
\BIBentryALTinterwordstretchfactor\fontdimen3\font minus
  \fontdimen4\font\relax}
\providecommand{\BIBforeignlanguage}[2]{{%
\expandafter\ifx\csname l@#1\endcsname\relax
\typeout{** WARNING: IEEEtran.bst: No hyphenation pattern has been}%
\typeout{** loaded for the language `#1'. Using the pattern for}%
\typeout{** the default language instead.}%
\else
\language=\csname l@#1\endcsname
\fi
#2}}
\providecommand{\BIBdecl}{\relax}
\BIBdecl

\bibitem{Doucette2011}
R.~T. Doucette and M.~D. McCulloch, ``A comparison of high-speed flywheels,
  batteries, and ultracapacitors on the bases of cost and fuel economy as the
  energy storage system in a fuel cell based hybrid electric vehicle,''
  \emph{Journal of Power Sources}, vol. 196, no.~3, pp. 1163--1170, 2011.

\bibitem{Guzzella2007}
L.~Guzzella, A.~Sciarretta \emph{et~al.}, \emph{Vehicle propulsion
  systems}.\hskip 1em plus 0.5em minus 0.4em\relax Springer, 2007, vol.~1.

\bibitem{Carrette2000}
L.~Carrette, K.~A. Friedrich, and U.~Stimming, ``Fuel cells: Principles, types,
  fuels, and applications,'' \emph{CHEMPHYSCHEM}, vol.~1, pp. 162--193, 2000.

\bibitem{Horiba2014}
T.~Horiba, ``Lithium-ion battery systems,'' \emph{proceedings of the IEEE},
  vol. 102, pp. 939--950, 2014.

\bibitem{Dhand2013}
A.~Dhand and K.~R. Pullen, ``Review of flywheel based internal combustion
  engine hybrid vehicles,'' \emph{International Journal of Automotive
  Technology}, vol.~14, no.~5, pp. 797--804, 2013.

\bibitem{Ju2014}
F.~Ju, Q.~Zhang, W.~Deng, and J.~Li, ``Review of structures and control of
  battery-supercapacitor hybrid energy storage system for electric vehicles,''
  in \emph{IEEE International Conference on Automation Science and
  Engineering}, Taipei, August 2014, pp. 18--22.

\bibitem{Shen2011}
C.~Shen, P.~Shan, and T.~Gao, ``A comprehensive overview of hybrid electric
  vehicles,'' \emph{International Journal of Vehicular Technology}, vol. 2011,
  2011.

\bibitem{Panday2014}
A.~Panday and H.~O. Bansal, ``A review of optimal energy management strategies
  for hybrid electric vehicle,'' \emph{International Journal of Vehicular
  Technology}, vol. 2014, 2014.

\bibitem{Paganelli2002}
G.~Paganelli, S.~Delprat, T.-M. Guerra, J.~Rimaux, and J.-J. Santin,
  ``Equivalent consumption minimization strategy for parallel hybrid
  powertrains,'' in \emph{Vehicular Technology Conference, 2002. VTC Spring
  2002. IEEE 55th}, vol.~4.\hskip 1em plus 0.5em minus 0.4em\relax IEEE, 2002,
  pp. 2076--2081.

\bibitem{Rousseau2008}
G.~Rousseau, D.~Sinoquet, A.~Sciarretta, and Y.~Milhau, ``Design optimisation
  and optimal control for hybrid vehicles,'' in \emph{EngOpt 2008 -
  International Conference on Engineering Optimization}, 2008.

\bibitem{Sciarretta2004}
A.~Sciarretta, M.~Back, and L.~Guzzella, ``Optimal control of parallel hybrid
  electric vehicles,'' \emph{IEEE Transactions on Control Systems Technology},
  vol.~12, no.~3, pp. 352--363, May 2004.

\bibitem{Baumann1998}
\BIBentryALTinterwordspacing
B.~Baumann, G.~Rizzoni, and G.~Washington, ``Intelligent control of hybrid
  vehicles using neural networks and fuzzy logic,'' in \emph{SAE Technical
  Paper}.\hskip 1em plus 0.5em minus 0.4em\relax SAE International, 02 1998.
  [Online]. Available: \url{http://dx.doi.org/10.4271/981061}
\BIBentrySTDinterwordspacing

\bibitem{Wang2006}
Z.~Wang, B.~Huang, W.~Li, and Y.~Xu, ``Particle swarm optimization for
  operational parameters of series hybrid electric vehicle,'' in \emph{Robotics
  and Biomimetics, 2006. ROBIO'06. IEEE International Conference on}.\hskip 1em
  plus 0.5em minus 0.4em\relax IEEE, 2006, pp. 682--688.

\bibitem{West2003}
M.~J. West, C.~M. Bingham, and N.~Schofield, ``Predictive control for energy
  management in all/more electric vehicles with multiple energy storage
  units,'' in \emph{Electric Machines and Drives Conference, 2003. IEMDC'03.
  IEEE International}, vol.~1, June 2003, pp. 222--228 vol.1.

\bibitem{Lee1998}
H.-D. Lee and S.-K. Sul, ``Fuzzy-logic-based torque control strategy for
  parallel-type hybrid electric vehicle,'' \emph{IEEE Transactions on
  Industrial Electronics}, vol.~45, no.~4, pp. 625--632, Aug 1998.

\bibitem{Lin2003}
C.-C. Lin, H.~Peng, J.~W. Grizzle, and J.-M. Kang, ``Power management strategy
  for a parallel hybrid electric truck,'' \emph{IEEE transactions on control
  systems technology}, vol.~11, no.~6, pp. 839--849, 2003.

\bibitem{Tate2000}
E.~D. Tate and S.~P. Boyd, ``Finding ultimate limits of performance for hybrid
  electric vehicles,'' SAE Technical Paper, Tech. Rep., 2000.

\bibitem{Perez2009}
\BIBentryALTinterwordspacing
L.~V. Pérez and E.~A. Pilotta, ``Optimal power split in a hybrid electric
  vehicle using direct transcription of an optimal control problem,''
  \emph{Mathematics and Computers in Simulation}, vol.~79, no.~6, pp. 1959 --
  1970, 2009, applied and Computational MathematicsSelected Papers of the Sixth
  PanAmerican WorkshopJuly 23-28, 2006, Huatulco-Oaxaca, Mexico. [Online].
  Available:
  \url{http://www.sciencedirect.com/science/article/pii/S0378475407001590}
\BIBentrySTDinterwordspacing

\bibitem{liu2008}
J.~Liu and H.~Peng, ``Modeling and control of a power-split hybrid vehicle,''
  \emph{Control Systems Technology, IEEE Transactions on}, vol.~16, no.~6, pp.
  1242--1251, Nov 2008.

\bibitem{Moura2011}
S.~Moura, H.~Fathy, D.~Callaway, and J.~Stein, ``A stochastic optimal control
  approach for power management in plug-in hybrid electric vehicles,''
  \emph{Control Systems Technology, IEEE Transactions on}, vol.~19, no.~3, pp.
  545--555, May 2011.

\bibitem{Kim2007}
M.-J. Kim and H.~Peng, ``Power management and design optimization of fuel
  cell/battery hybrid vehicles,'' \emph{Journal of Power Sources}, vol. 165,
  pp. 819--832, 2007.

\bibitem{Piccolo2001}
A.~Piccolo, L.~Ippolito, V.~zo~Galdi, and A.~Vaccaro, ``Optimisation of energy
  flow management in hybrid electric vehicles via genetic algorithms,'' in
  \emph{Advanced Intelligent Mechatronics, 2001. Proceedings. 2001 IEEE/ASME
  International Conference on}, vol.~1.\hskip 1em plus 0.5em minus 0.4em\relax
  IEEE, 2001, pp. 434--439.

\bibitem{Qi2016}
X.~Qi, G.~Wu, K.~Boriboonsomsin, and M.~J. Barth, ``Development and evaluation
  of an evolutionary algorithm-based online energy management system for
  plug-in hybrid electric vehicles,'' \emph{IEEE Transactions on Intelligent
  Transportation Systems}, vol.~PP, no.~99, pp. 1--11, 2016.

\bibitem{Qi2016data}
X.~Qi, G.~Wu, K.~Boriboonsomsin, M.~J. Barth, and J.~Gonder, ``Data-driven
  reinforcement learning--based real-time energy management system for plug-in
  hybrid electric vehicles,'' \emph{Transportation Research Record: Journal of
  the Transportation Research Board}, no. 2572, pp. 1--8, 2016.

\bibitem{Musardo2005}
C.~Musardo, G.~Rizzoni, Y.~Guezennec, and B.~Staccia, ``A-ecms: An adaptive
  algorithm for hybrid electric vehicle energy management,'' \emph{European
  Journal of Control}, vol.~11, no. 4-5, pp. 509--524, 2005.

\bibitem{Roy2014}
H.~Roy, A.~McGordon, and P.~Jennings, ``A generalized powertrain design
  optimization methodology to reduce fuel economy variability in hybrid
  electric vehicles,'' \emph{Vehicular Technology, IEEE Transactions on},
  vol.~63, no.~3, pp. 1055--1070, March 2014.

\bibitem{Kaiser1944}
F.~Kaiser, ``The climb of jet-propelled aircraft, part i. speed along the path
  of optimum climb,'' Ministry of Supply (Gt. Brit.), Tech. Rep., April 1944.

\bibitem{Tsien1951}
H.~S. Tsien and R.~C. Evans, ``Optimal thrust programming for a sounding
  rocket,'' \emph{Journal of the American Rocket Society}, vol.~21, no.~3, pp.
  99--107, 1951.

\bibitem{Perantoni2015Track}
G.~Perantoni and D.~Limebeer, ``Optimal control of a formula one car on a
  three-dimensional track. {P}art 1: Track modelling and identification,''
  \emph{Journal of Dynamic Systems, Measurement, and Control}, vol. 137, pp.
  051\,018--051\,018--11, May 2015.

\bibitem{Pacejka2006}
H.~B. Pacejka, \emph{Tyre and Vehicle Dynamics, 2nd Edition}.\hskip 1em plus
  0.5em minus 0.4em\relax Butterworth-Heinemann Ltd, 2006.

\bibitem{Limebeer20153DOptimal}
D.~Limebeer and G.~Perantoni, ``Optimal control of a formula one car on a
  three-dimensional track�part 2: Optimal control,'' \emph{Journal of Dynamic
  Systems, Measurement, and Control}, vol. 137, pp. 051\,019--051\,019--13,
  2015.

\bibitem{ADVISOR2}
N.~R.~E. Laboratory, ``Advanced vechile simulator,''
  http://adv-vehicle-sim.sourceforge.net/advisor\_doc.html, 2003.

\bibitem{Doucette2013}
R.~Doucette, ``The oxford vehicle model: A tool for modeling and simulating the
  powertrains of electric and hybrid electric vehicles,'' Ph.D. dissertation,
  University of Oxford, 2013.

\bibitem{Patterson2013}
M.~A. Patterson and A.~V. Rao, \emph{GPOPS - II Version 1.0: A General-Purpose
  MATLAB Toolbox for Solving Optimal Control Problems Using the Radau
  Pseudospectral Method}, University of Florida, 2013.

\bibitem{Abramowitz1965}
M.~Abramowitz and I.~Stegun, \emph{Handbook of Mathematical Functions with
  Formulas, Graphs, and Mathematical Tables}.\hskip 1em plus 0.5em minus
  0.4em\relax Dover Publications, 1965.

\bibitem{Darby2011}
C.~L. Darby, W.~W. Hager, and A.~V. Rao, ``An hp-adaptive pseudospectral method
  for solving optimal control problems,'' \emph{OPTIMAL CONTROL APPLICATIONS
  AND METHODS}, vol.~32, pp. 476--502, 2011.

\bibitem{Biegler2008}
L.~Biegler and V.~Zavala, ``Large-scale nonlinear programming using ipopt: An
  integrating framework for enterprise-wide dynamic optimization,''
  \emph{Computers and Chemical Engineering}, vol.~33, pp. 575--582, 2008.

\bibitem{Weinstein2014}
\BIBentryALTinterwordspacing
M.~J. Weinstein and A.~V. Rao, ``A source transformation via operator
  overloading method for the automatic differentiation of mathematical
  functions in matlab,'' \emph{ACM Transactions on Mathematical Software}, vol.
  In Revision, 2014. [Online]. Available:
  \url{http://vdol.mae.ufl.edu/SubmittedJournalPublications/gatorad-TOMS-November-2013.pdf}
\BIBentrySTDinterwordspacing

\bibitem{Rao2010}
A.~V. Rao, D.~A. Benson, D.~A. Benson, C.~Darby, M.~A. Patterson, C.~Francolin,
  I.~Sanders, and G.~T. Huntington, ``Algorithm 902: Gpops, a matlab software
  for solving multiple-phase optimal control problems using the gauss
  pseudospectral method,'' vol.~37, no.~22, 2010.

\end{thebibliography}

\appendix
\begin{table}[htb]
\caption{Vehicle parameter values}
\label{vehicleParameters}
\begin{center}
\small
\begin{tabular}{|c||c||c|}
\hline
Symbol & Description & Value \\ 
\hline
$M$ & Vehicle total mass & 1400 $kg$ \\ 
$M_{car}$ & Vehicle only mass & 1280 $kg$ \\ 
$M_{bat}$ & Battery and motor mass & 100 $kg$ \\ 
$M_{fly}$ & Flywheel mass & 20 $kg$ \\ 
$I_x$ & x moment of inertia & 500 $kg \, m^2$\\
$I_y$ & y moment of inertia & 1000 $kg \, m^2$\\
$I_z$ & z moment of inertia & 1000 $kg \, m^2$\\
$a$ & Mass centre from front axle &1.35 $m$ \\
$b$ & Mass centre from rear axle & 1.35 $m$  \\
$h$ & Centre of mass height & 0.5 $m$ \\
$C_d$ & Drag coefficient & 0.3 \\
$A$ & Vehicle frontal area & 1.8 $m^2$ \\
$\rho$ & Air density & 1.2 $kg/m^3$ \\
$C_d$ & Rolling resistance coefficient & 0.009 \\
$E_{b}^{max}$ & Max battery energy capacity & 5 $MJ$ \\
$P_{b}^{max}$ & Max battery power & 25 $kW$ \\
$V_{OC}^{min}$ & Min battery voltage & 240 $V$ \\
$V_{OC}^{max}$ & Max battery voltage & 210 $V$ \\
$R_b$ & Battery internal resistance & 0.5 $\Omega$ \\
$SoC_{min}$ & Min state of charge & 40\,\% $V$ \\
$SoC_{max}$ & Max state of charge & 80\,\% $V$ \\
$P_{fly}^{max}$ & Max flywheel power & 60 $kW$ \\
$\mu_{em}$ & Electric motor efficiency & 85\% \\
$E_{fly}^{max}$ & Max flywheel energy & 400 $kJ$ \\
$J_{fly}$ & Flywheel spinning inertia & 0.02 $kg \, m^2$ \\
$\mu_{CVT}$ & CVT efficiency & 85\% \\

\hline 
\end{tabular}
\end{center}
\end{table}

\begin{IEEEbiography}[{\includegraphics[width=1in]{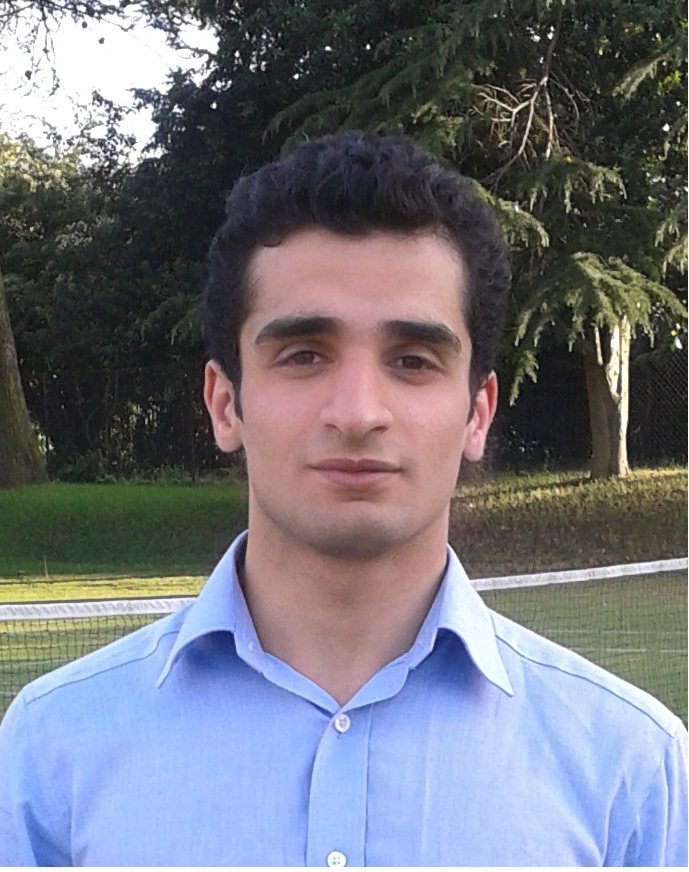}}]{Mehdi Imani Masouleh}
received the M.Eng. degree in electrical and electronic engineering with first-class honours from the
Imperial College London, London, U.K., in 2011. He is currently pursuing the D.Phil. degree in
engineering science at the University of Oxford.
His current research interests include optimal control, vehicle dynamics, nonlinear stability and hybrid vehicles.
\end{IEEEbiography}

\begin{IEEEbiography}[{\includegraphics[width=1in]{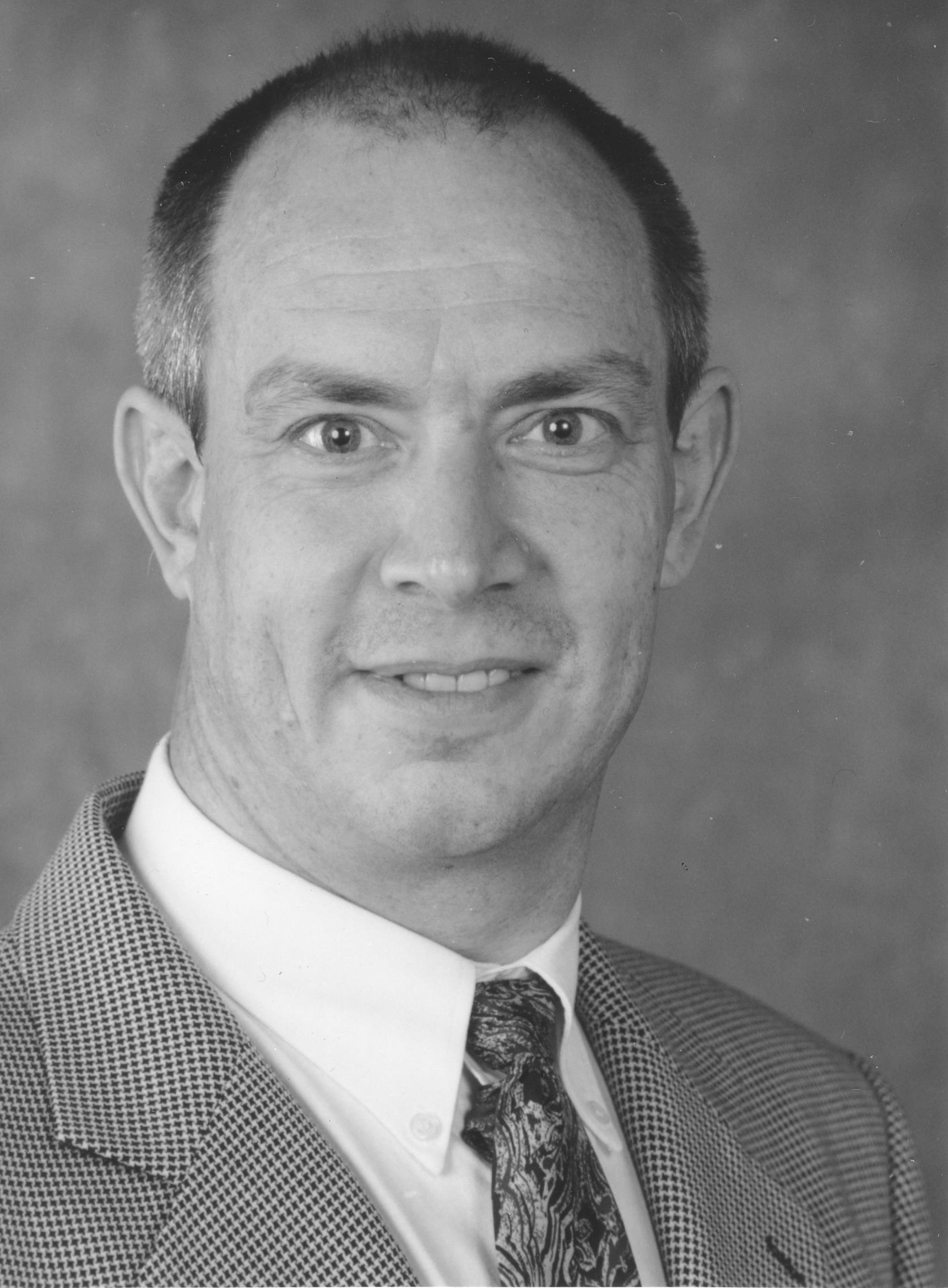}}]{David J N Limebeer}
received a B.Sc.(Eng) degree from the University of the Witwatersrand in 1974, MSc(Eng) and PhD degrees from the University of Natal in 1977 and 1980, respectively, and the DSc (Eng) from the University of London in 1992. He was a post-doc researcher at the University of Cambridge between 1980 and 1984. He then joined the Electrical and Electronic Engineering Department at Imperial College as a lecturer. He was promoted to Reader in 1989, Professor in 1993, Head of the Control Group in 1996, and Head of Department 1999-2009. In 2009 he moved to Oxford as Professor of Control Engineering and Professorial Fellow at New College Oxford. His research interests include applied and theoretical problems in control systems and engineering dynamics. He is a Fellow of the IEEE (1992), a Fellow of the IET (1994), and a Fellow of the Royal Academy of Engineering (1997).
\end{IEEEbiography}
\end{document}